\documentclass[letterpaper,11pt,reqno]{amsart}
\usepackage{amsmath}
\usepackage{amssymb}
\usepackage{amsfonts}
\usepackage{amssymb}	

\usepackage{bm}

\usepackage{amssymb}
\usepackage{amsmath}
\usepackage{amsthm}

\usepackage{graphicx}

\usepackage{enumerate}	

\usepackage{tikz}				
\usetikzlibrary{arrows,matrix,positioning}
\usetikzlibrary{calc}

\newtheorem{theorem}{Theorem}[section] 
\newtheorem{lemma}[theorem]{Lemma}
\newtheorem{corollary}[theorem]{Corollary}

\theoremstyle{definition}

\theoremstyle{definition}

\theoremstyle{plain}

\theoremstyle{definition}
\newtheorem{remark}[theorem]{Remark}

\theoremstyle{definition}

\theoremstyle{definition}

\numberwithin{figure}{section}
\numberwithin{theorem}{section}

\newcommand{\ColorMap}{\mathcal{C}}

\newcommand{\C}{\mathbb{C}}

\newcommand{\R}{\mathbb{R}}

\renewcommand{\H}{\mathbb{H}}

\newcommand{\F}{\mathcal{F}}

\newcommand{\D}{\mathbb{D}}  

\newcommand{\EEL}{\mathbb{E}^{+}}

\newcommand{\Esub}{\hbox{\tiny$\EEL$}}

\newcommand{\MG}{\mathsf{S}_{\H}}  

\newcommand{\MGZ}{\Gamma}  

\newcommand{\Z}{\mathbb{Z}}

\newcommand{\E}{\mathsf{E}}

\newcommand{\T}{T}

\newcommand{\I}{\mathcal{I}}

\newcommand{\M}{\mathcal{M}}

\renewcommand{\Im}[1]{\mathrm{Im}\, #1}

\renewcommand{\Re}[1]{\mathrm{Re}\, #1}
\newcommand{\Sur}[2]{S^{#1}_{#2}}
\newcommand{\Ru}[1]{\mathcal{R}_{#1}}

\newcommand{\Xrightarrow}[1]{\xrightarrow{\ \ \; #1\ \ }}

\newcommand{\Mod}[1]{\ \, \mathrm{mod}\ #1}

\renewcommand{\S}{\mathsf{S}}

\newcommand{\Fig}{Figure}

\newcommand{\Emph}[1]{\emph{\textbf{#1}}}

\title{Symmetry and Art}

\author{Emily J.~Gullerud}
\author{James S.~Walker}

\address{Emily J.~Gullerud\\
         School of Mathematics\\
         Univ.~of Minnesota\\
         gulle069@umn.edu}

\address{James S.~Walker\\
         Department of Mathematics\\
         Univ.~of Wisconsin-Eau Claire\\
         walkerjs@uwec.edu}

\begin{document}

\vspace*{-2.3cm}
\maketitle

\begin{abstract}
We use some fundamental ideas from complex analysis to create symmetric images and animations. Using a domain coloring algorithm, we generate mappings to the entire complex plane or the hyperbolic upper half-plane. The resulting designs can have rotational, translational, or mirror symmetry according to our chosen mapping functions. An appealing feature of these designs is
how they reveal important properties of Euclidean and non-Euclidean geometries. We can also generate animations of our designs. Our goal is to create designs and animations having significant artistic content.
\end{abstract}

\begin{quotation}\textsf{\textsl{\small
\dots there's mathematical beauty in a piece that's at least as important as the
surface beauty that everyone can see.}}---Frank Farris
\end{quotation}

\section{Introduction}\thispagestyle{empty}
The connection between symmetry and artistic design has a long history (see~e.g.,Weyl~\cite{ref:Weyl}).
In this paper we describe the application of domain coloring to the creation of symmetric designs and animations. First, we examine domain coloring in its original setting as a method for sketching graphs of complex-valued functions. Then we
will describe the creation of designs in the Euclidean complex plane and the non-Euclidean hyperbolic upper half-plane.

\begin{figure}[!b]\centering
\setlength{\unitlength}{1in}
\begin{picture}(2.5,2.5)
\put(0,0){\resizebox{!}{2.5in}{
\begin{picture}(3.03,3.03)
\put(0,1.5){\includegraphics[width=1.48in]{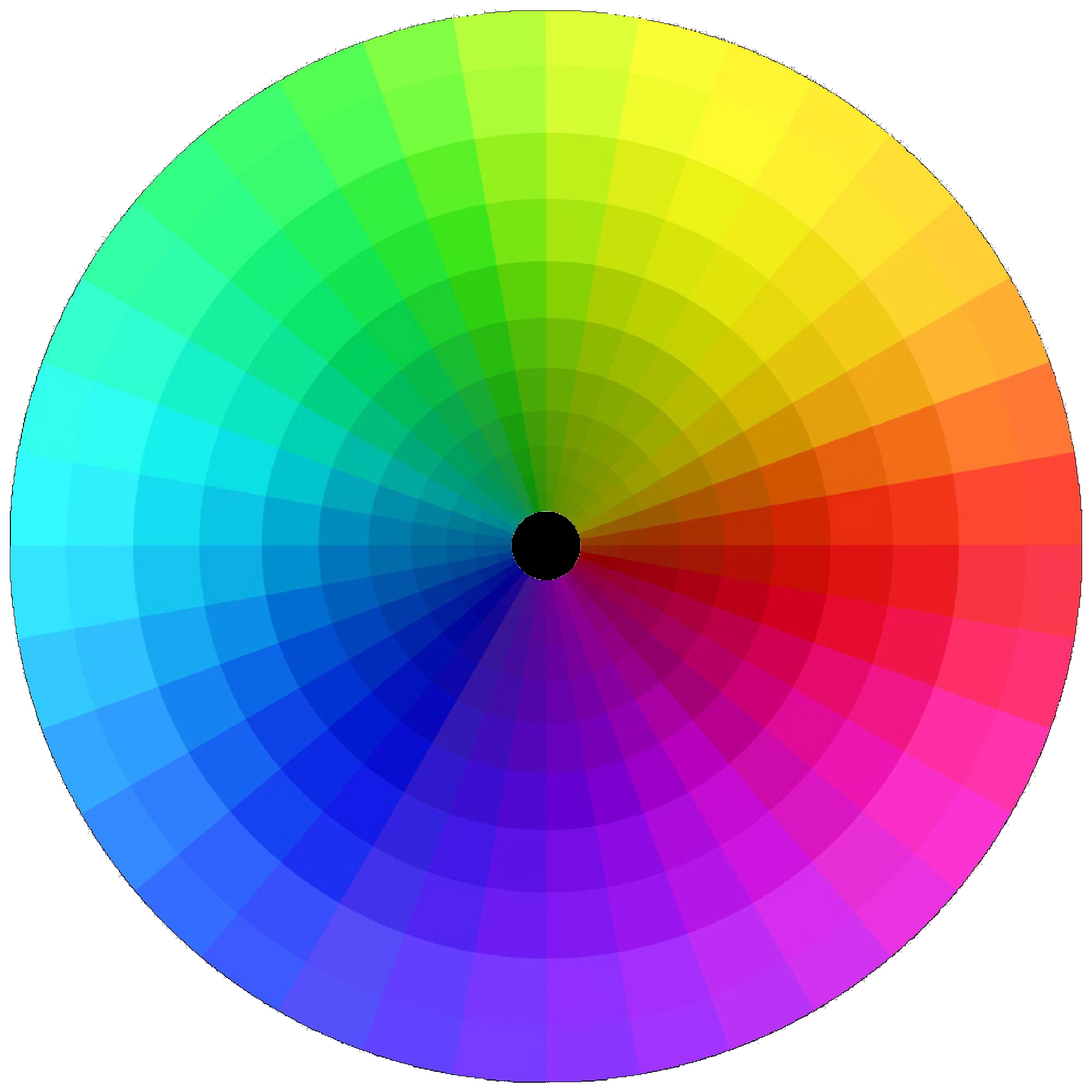}}
\put(0,1.5){\framebox(1.48,1.48){}}
\put(1.55,1.5){\includegraphics[width=1.48in]{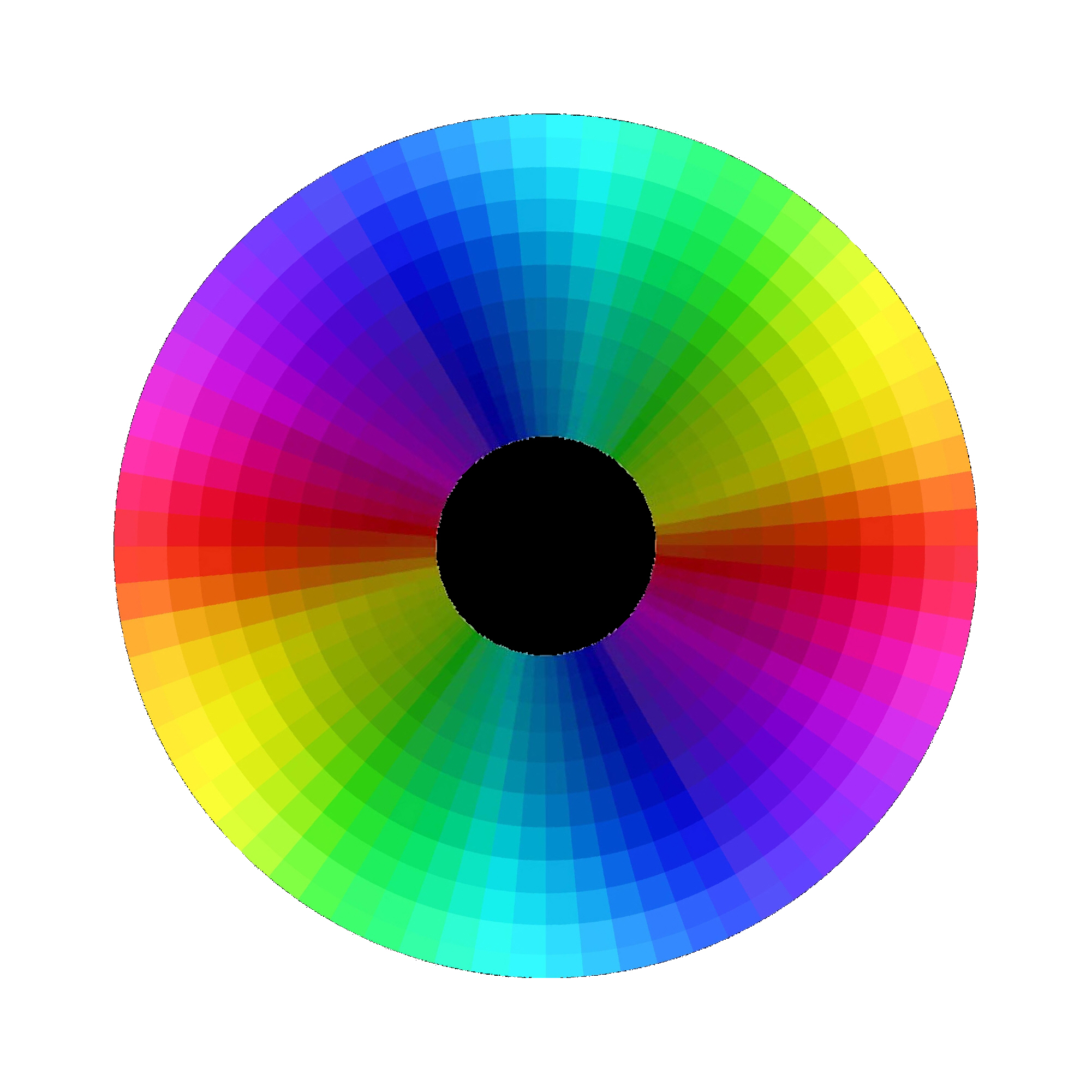}}
\put(1.55,1.5){\framebox(1.48,1.48){}}
\put(0,0){\includegraphics[width=3.03in]{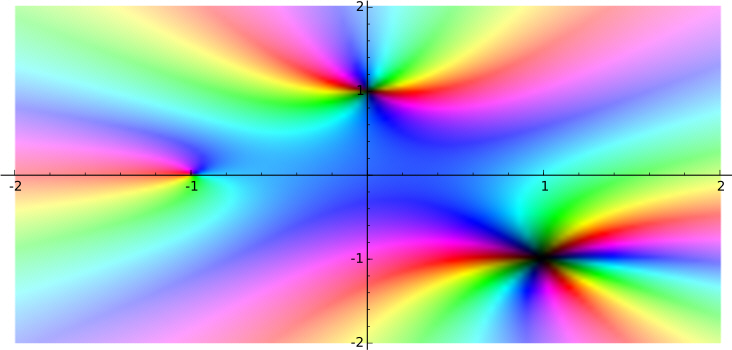}}
\end{picture}}}
\end{picture}
\caption{\small Top left: Color plot for $w = z$ over $[-2, 2] \times [-2,2]$; a basic color map $\ColorMap$.
Top right: Color plot for $w = z^2$. Bottom: SageMath color plot for $w = 3 (z+1) (z-i)^2 (z-1+i)^3$.
{\scriptsize First published in ``Creating Symmetry: The Artful Mathematics of Wallpaper Patterns'' by James S. Walker, Notices of the Amer.~Math.~Soc., Vol~62, No~11 (Dec.~2015), pp.~1350--1354, published by the American Mathematical Society.
\copyright 2015 American Mathematical Society. Used by permission.} 
\label{fig:ColorWheelPlots}}
\end{figure}
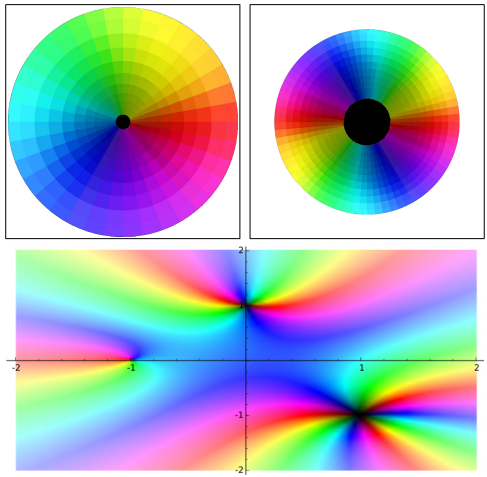

At the very end of the twentieth century, Farris~\cite{ref:FarrisColorWheel} played a principal role in developing the use of domain coloring to sketch graphs of functions, $f\colon\C\to\C$. The idea is to use a \emph{color wheel}, a well-known tool in the visual arts. A basic example of a color wheel is shown at the top left of \Fig~\ref{fig:ColorWheelPlots}. This color wheel is used to mark locations in the complex plane. For each value of $w$, the value of $\ColorMap(w)$ is a unique color (at least in principle). For example, on the top left of \Fig~\ref{fig:ColorWheelPlots}, the values of $w$ that are near $i$ are colored greenish-yellow, while values near $-1$ have a light blue tint. As values approach zero they turn black, and beyond a certain radius they are all colored white. By composing a function $f(z)$ with $\ColorMap$, we get a function $\ColorMap(f(z))$ that gives a color portrait of $f(z)$. For instance, on the top right of \Fig~\ref{fig:ColorWheelPlots}, a color map for $w = z^2$ is shown. Notice how the colors near $w = 0$ cycle twice through the rainbow as we move once around $w = 0$. The connection to winding numbers is thereby made visually evident. Also, contour lines meet at right angles away from $w = 0$, just as they do for the function $w = z$. This illustrates conformality of $w = z^2$ away from the origin.
There is software now that produces color plots with great ease~\cite{ref:SageMath, ref:SymmetryWorks}. At the bottom of \Fig~\ref{fig:ColorWheelPlots}, we show a color plot of the function $w = 3(z+1) (z-i)^2 (z-1+i)^3$. This plot was produced with the free SageMath system~\cite{ref:SageMath,ref:SMComplexPlot}. Just these two commands were needed:
\[\hbox{
\begin{tabular}{l}
\verb!3*(z+1)*(z-i)^2*(z-1+i)^3!\\
\verb!complex_plot(f, (-2, 2), (-2, 2))!\\
\end{tabular}}
\]
The plot that SageMath produced clearly marks the location of the zeros at $-1$, $i$, and $1 - i$ and their multiplicities of $1$, $2$, and $3$, respectively. All this information is encoded in the number of times the colors of the rainbow are cycled through in the neighborhood of each zero. There are several nice examples of color plots at the web~site created by Crone~\cite{ref:CroneDomainColoring}, including plots of branching in Riemann surfaces. The reader may also wish to
explore color plotting using the free software, \textsc{SymmetryWorks}~\cite{ref:SymmetryWorks}.
\subsection{Pseudocode for Creating Color Plots\label{subsec:PseudoCode}}
\begin{figure}[!t]\centering
\setlength{\unitlength}{1in}
\begin{picture}(4.65,1.2)
\put(0,0){\includegraphics[width=1.45in]{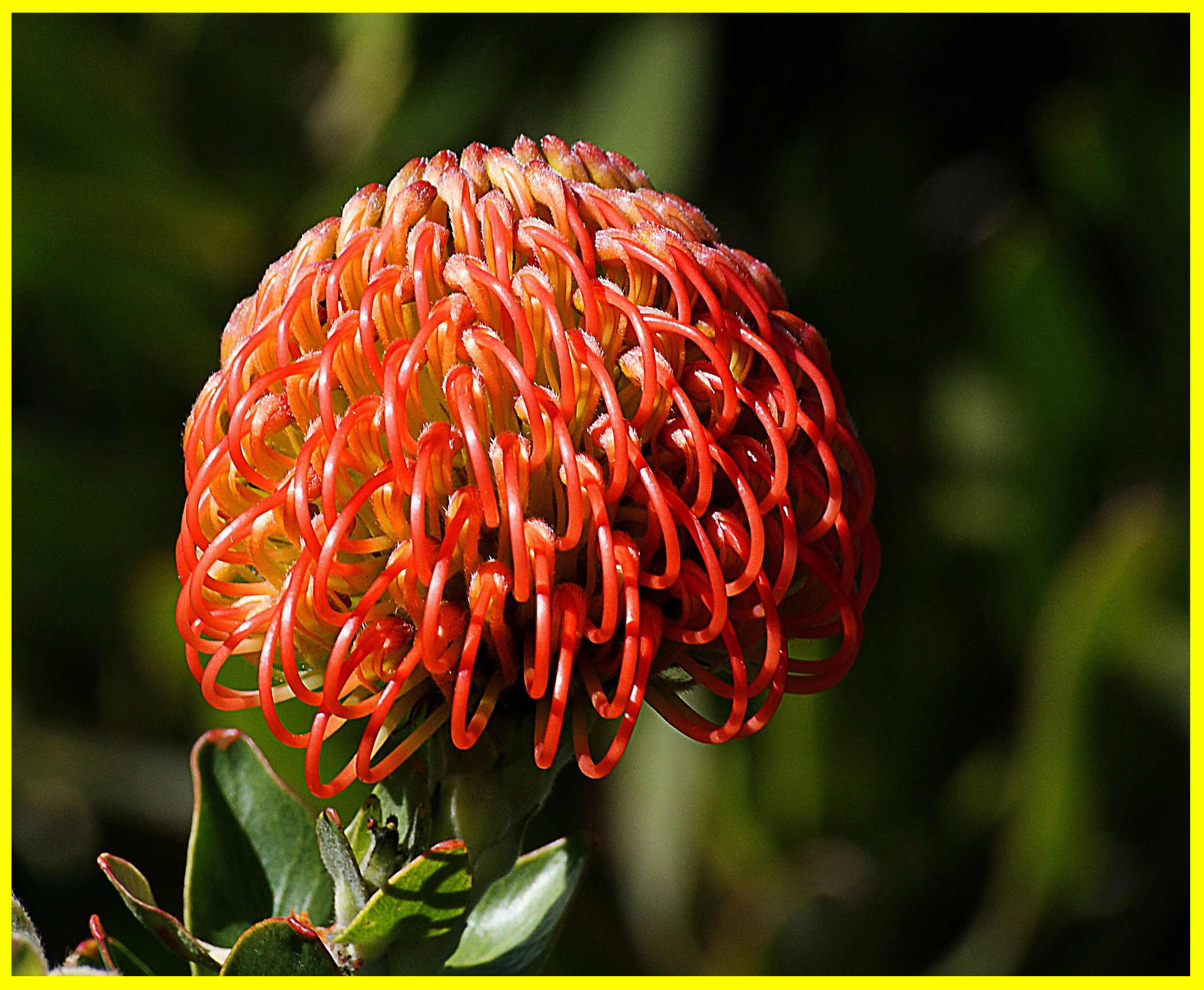}}
\put(1.6,0){\includegraphics[width=1.45in]{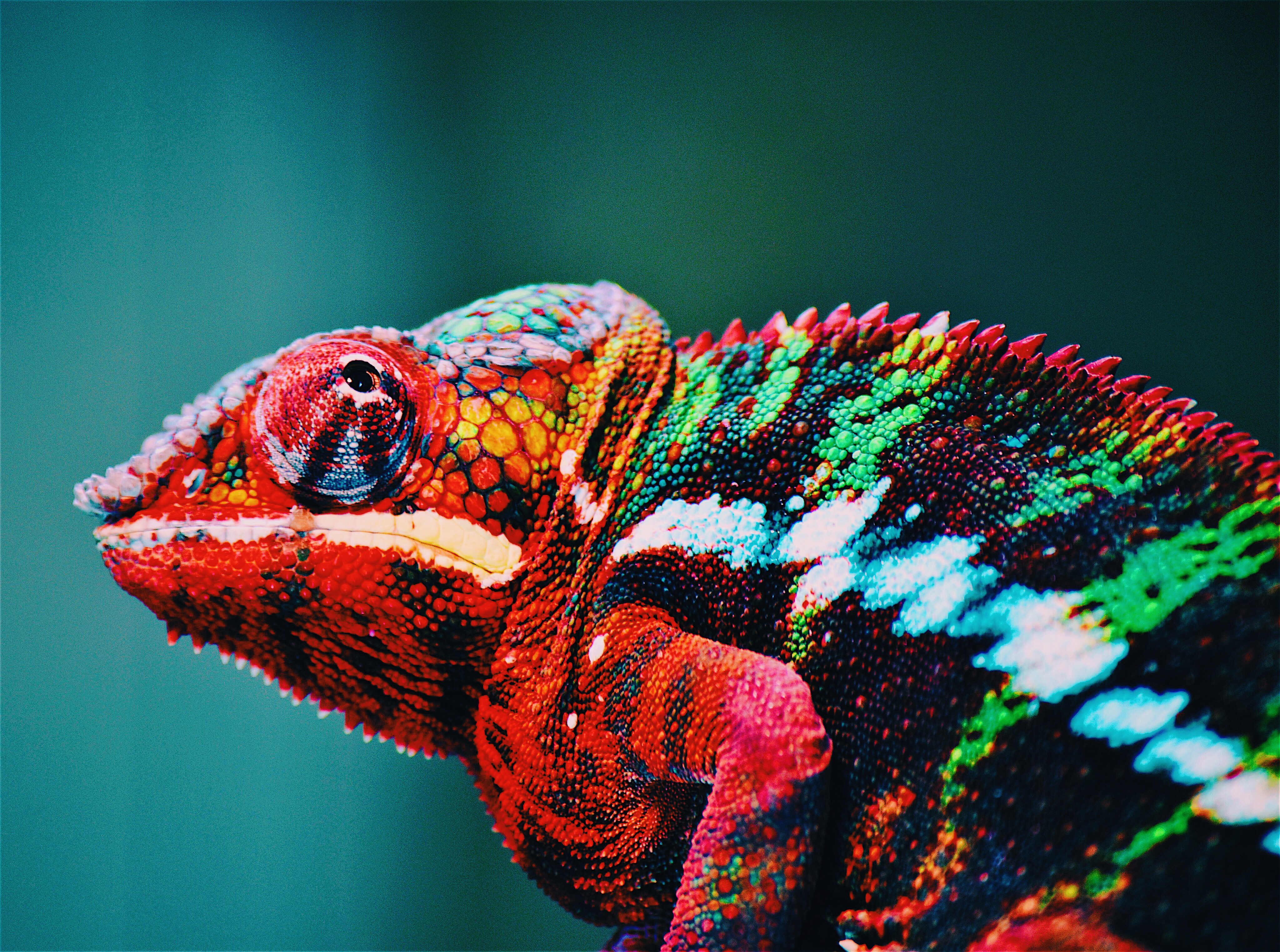}}
\put(3.2,0){\includegraphics[width=1.45in]{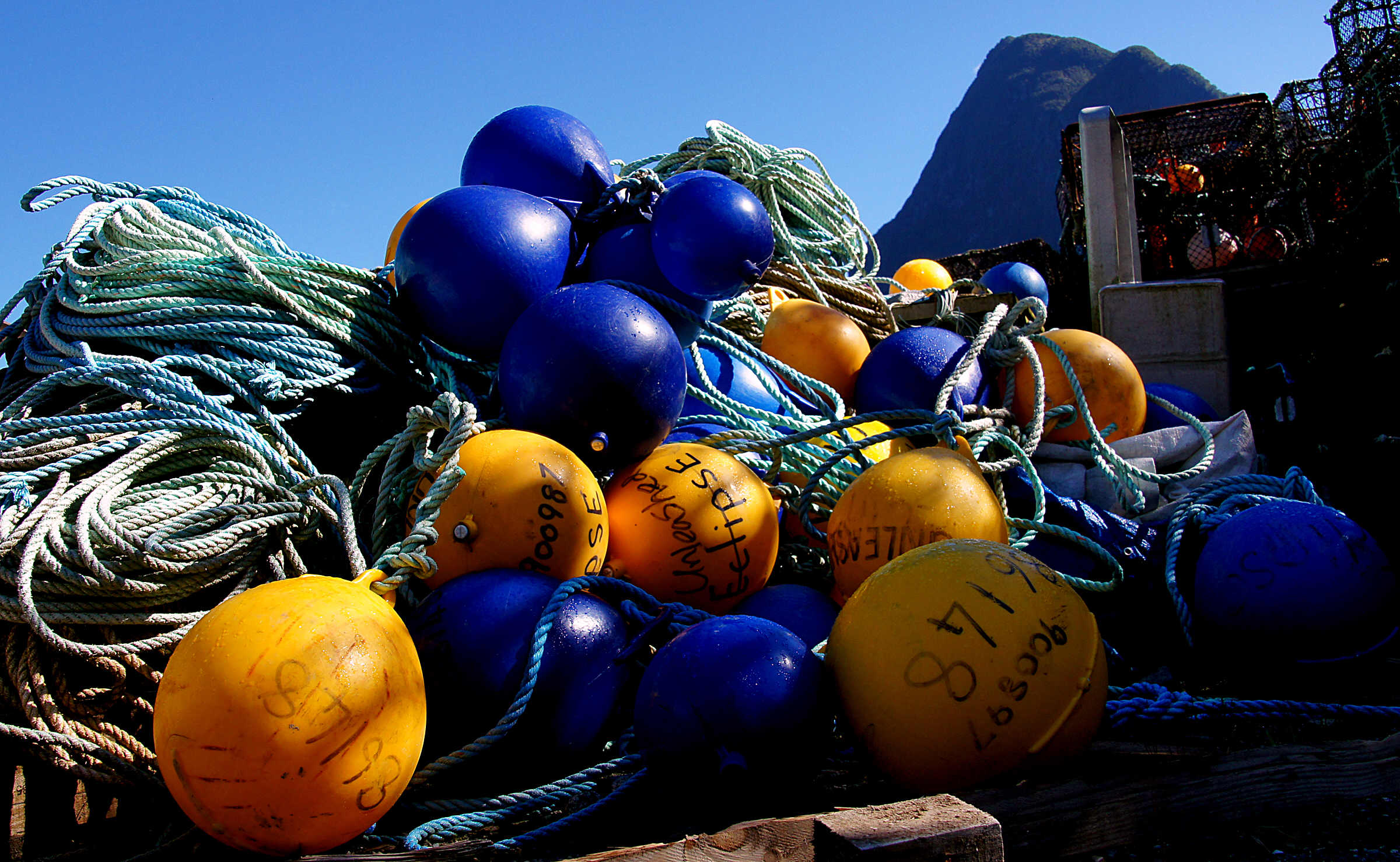}}
\end{picture}
\caption{\small Images for color maps, from public domain stock images found on the Internet.
Left to right: \emph{Waratah flower} with yellow border, \emph{Reptile}, \emph{Buoy}. 
\label{fig:SourceImages}}
\end{figure}

We have found the most challenging part of producing color plots is the
creation of a domain colored image using a specific color map. In Figure~\ref{fig:ColorMapProcedure}
we provide pseudocode for this procedure. It is usually a straightforward
task to supplement this procedure with additional code for creating symmetric
designs using various functions.
\begin{figure}[!htb]\centering
{\small\textbf{Domain Coloring Procedure}}
\begin{footnotesize}
\begin{verbatim}
CMp = loadimage(filename) %load the RGB image color map
R = number of rows of CMp and C = number of columns of CMp
ScaleFactor = 10 %a scaling factor for units in color map
mCent = R / 2 and nCent = C / 2 %coordinates of image center
xinc = (2*ScaleFactor) / R and yinc = (2*ScaleFactor) / C %x, y increments
%Create color-mapped display.
M = 1000 %increase M for higher resolution
L = 2  %increase/decrease L to zoom out/in
incr = 2*L/M  and  [x, y] = rectangular grid of points from -L to L spaced by incr
Initialize wImg %Matrix wImg, all zeros for RGB values, M rows, M columns.
z = x + iy  %create array of complex values
w = f(z) %compute values of the function f over grid
s = real(w) and t = imag(w) %real and imaginary parts of function values
for m = 1 to M
  for n = 1 to M
    %Calculate coordinates  p and q for color map
    p = mCent + round(s(m,n)/xinc)  and  q = nCent + round(t(m,n)/yinc)
    if p>0 and p<R then
      if q>0 and q<C then
      wImg(m,n) = CMp(p,q) %assign RGB values from color map
      end if
    end if
  end for
end for
Display(wImg) %Display RGB image wImg
\end{verbatim}
\end{footnotesize}
\caption{\small Pseudocode for applying a color map to function $w = f(z)$.\label{fig:ColorMapProcedure}}
\end{figure}
For color maps, we will adopt the clever idea of Farris of using color photos of natural scenes, rather than
a simple color wheel. The color images we used for the designs in this paper are shown in
Figure~\ref{fig:SourceImages}. For those readers not wishing to reinvent the wheel, the software~\cite{ref:SymmetryWorks}
can be used for artistic designs like the ones we discuss in the remainder of this paper.

\section{Designs in the Euclidean Complex Plane}
In the geometry of the Euclidean plane, there are three basic symmetry operations (congruences). They
are (1) translation, (2) rotation, (3) reflection. For a complex variable, $z = x + iy$, these symmetry operations
are exemplified by the following mappings:
\begin{enumerate}
\item Translation:\ \ $\T_c:\ z\longrightarrow z + c$, where $c = a + i b$ is a complex constant.

\item Reflection through the $x$-axis:\ \ $R_x:\ z\longrightarrow \overline{z}$, i.e., $z = x+iy \longrightarrow x - iy$.

\item Rotation by $\theta$ around the origin:\ \ $\rho_\theta:\ z \longrightarrow e^{i\theta}\, z$, where $e^{i\theta} = \cos \theta + i\sin \theta$.
\end{enumerate}
Any Euclidean symmetry (congruence) can be expressed as a finite composition of these basic mappings. This set of
all symmetries of the plane is a mathematical group, the \Emph{Euclidean group}, denoted by $\E_2$.

For creating symmetric designs, we will apply domain coloring to a function $w = f(z)$ that has symmetry with respect
to a group of transformations. If $S$ is a transformation of the complex plane $\C$, then the function $f$ is symmetric
with respect to $S$ if
\begin{equation}\label{eq:SymmetryOfFunction}
f\bigl(S(z)\bigr) = f(z)\quad\text{for all $z\in\C$}
\end{equation}
The set of symmetries of a function $f$ is a mathematical group, $\S_f$. Our method for creating designs is to apply
domain coloring to functions that have some pre-assigned symmetry group. For that purpose,
the Euclidean group is too large. The only functions symmetric under all mappings in
$\E_2$ are constant functions. To obtain non-trivial symmetric designs we restrict to discrete subgroups of $\E_2$.
We will examine symmetric designs created from symmetries with respect to these three discrete subgroups:
\begin{enumerate}
\item Symmetry with respect to a lattice of two-dimensional translations
\[
\T_{mu + nv}: z \longrightarrow z + mu + nv\quad \text{all $m, n \in \Z$},
\]
where $u$ and $v$ are non-collinear complex numbers. Enforcing this symmetry creates a tessellation over the
lattice in $\C$ generated by $u$ and $v$. 
\item Symmetry with respect to rotation about the origin by an $n$th root of unity:
\[
\rho_{2\pi/n}: z \longrightarrow e^{i2\pi/n} z\quad \text{for fixed $n\in\Z, n \neq 0$}
\]
Enforcing this symmetry creates a \Emph{rosette}, having $n$-fold rotational symmetry about
the origin. 
\item Symmetry with respect to both a lattice of two-dimensional translations and rotation by an $n$th root of
unity about the origin. Symmetric designs of this type are only possible for
rotations $\rho_{2\pi/n}$ with $n = 2$, $3$, $4$, or $6$. We shall refer to a design with this symmetry
as a \Emph{rotationally symmetric wallpaper pattern}.
\end{enumerate}
All of these basic examples can be modified to add more symmetry to a design, including various reflective symmetries. We will discuss these modifications
as we examine these three basic symmetric designs.
\subsection{Translational Symmetry\label{subsec:LatticeBasedSymmetry}}
\begin{figure}[!htb]\centering
\setlength{\unitlength}{1in}
\begin{picture}(4.84,2.25)
\put(0,0){\resizebox{!}{2.2in}{
\begin{picture}(5.8,2.5) 
\put(0,0){\resizebox{3in}{2.5in}{\includegraphics{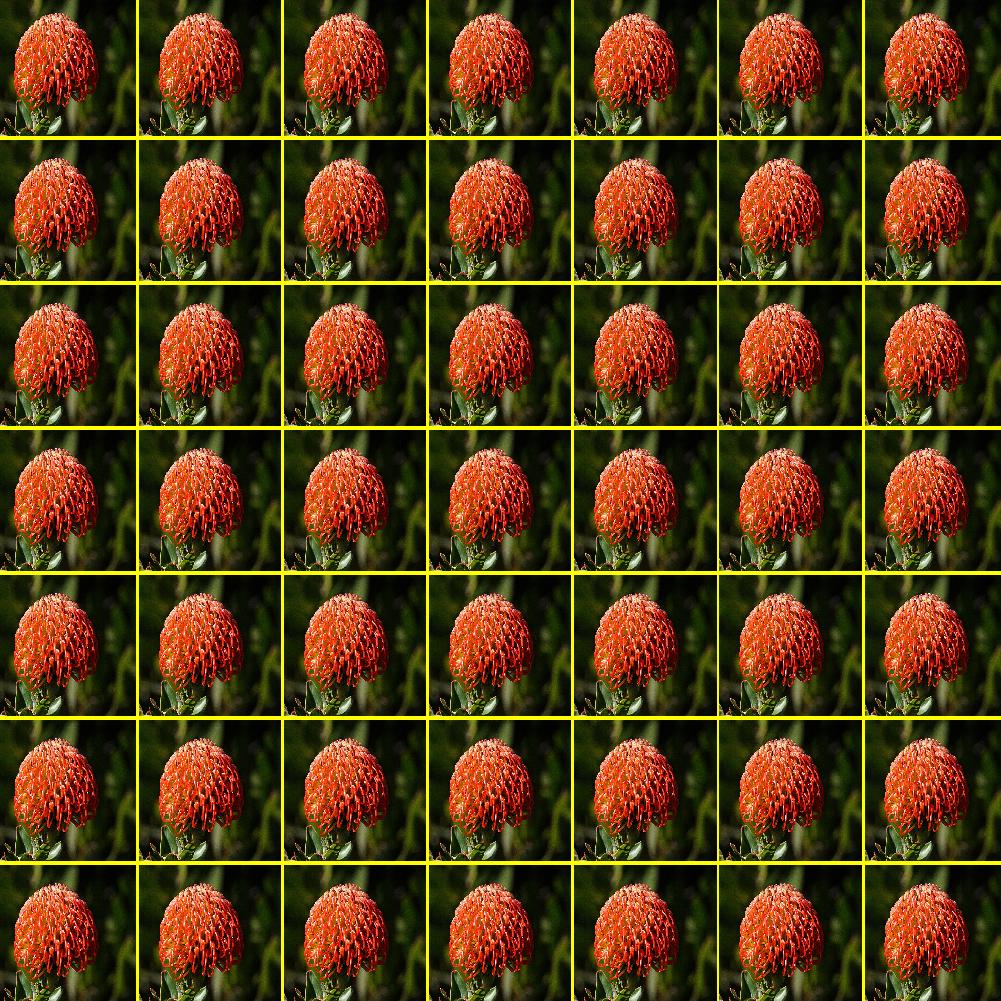}}}
\put(3.15,0){\resizebox{!}{2.5in}{\includegraphics{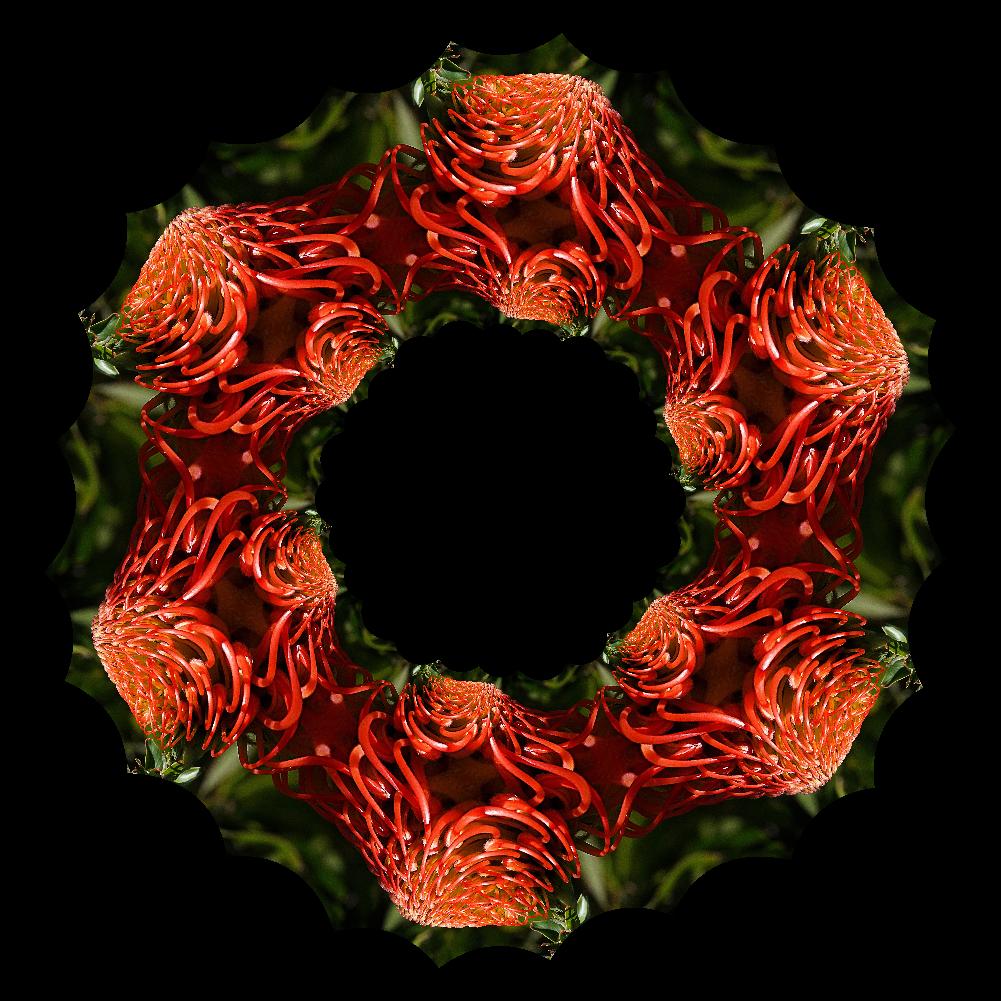}}}
\end{picture}}}
\end{picture}
\caption{\small Left: Design having translational symmetry over a square lattice.
Right: Rosette design having 6-fold rotational symmetry.\label{fig:TesselationAndRosetteDesigns}}
\end{figure}

On the left of \Fig~\ref{fig:TesselationAndRosetteDesigns}, we show
a design having symmetry with respect to a lattice of two-dimensional translations. The lattice is
rectangular, and is generated by the non-collinear complex numbers $u = 2100$ and $v = 1700i$. 
The positive integers, $2100$ and $1700$,  
are the width and height in pixels of the image used for
the color map. That color map image is the \emph{Waratah flower} image shown in \Fig~\ref{fig:SourceImages}.
The function $f(z)$ that we used to creates the design is, for $z = x + i y$,
\begin{equation}\label{eq:SquareLatticeSymmetricDesign}
f(x + i y) = [x\Mod 2100] +  i [y\Mod 1700]
\end{equation}
This formula assumes that the lower left corner of the \emph{Waratah flower} image is located at the origin
in $\C$, and uses mod-equivalence from elementary number theory to ensure that $f$ is symmetric with
respect to translations over the lattice generated by $u$ and $v$.
\noindent%
The image on the left of \Fig~\ref{fig:TesselationAndRosetteDesigns} resembles a sheet of postage stamps, rather than a complex artistic design. Although if one began with an image of a person (say Marilyn Monroe), or a commodity (say a Campbell's soup can) then one could easily produce designs reminiscent of some of Andy Warhol's famous pictures.
\subsection{Rotational Symmetry\label{subsec:RotationallySymmetricDesign}}
\subsubsection{Rosettes\label{subsubsec:Rosettes}}
On the right of \Fig~\ref{fig:TesselationAndRosetteDesigns}, we show an example of a design having
symmetry with respect to rotation about the origin by a $6$th root of unity. In other words, it has
$6$-fold rotational symmetry. The color map used to create the image is the \emph{Waratah flower} image
shown in \Fig~\ref{fig:SourceImages}, \emph{but without the yellow border}. The function we used to
create the design is
\begin{equation}\label{eq:RosetteDesignFunction}
f(z) = (1+i) + i\, \frac{z^6}{4}  + z^{-6}
\end{equation}
Since $\left(e^{i2\pi / 6}\right)^6 = 1$, it follows that $f(z)$ is symmetric with respect to the rotation
$\rho_{2\pi/6}$, i.e., it has $6$-fold rotational symmetry. We chose array values for $z = x + iy$ so that
the rosette design is a square image, of dimensions $1000 \times 1000$ pixels. 
The black colored regions in the design correspond to values of $f$ lying outside the region in $\C$ corresponding to the color map. Later, we will show an interesting alternative way to color these values.

\begin{figure}[!htb]\centering
\setlength{\unitlength}{1in}
\begin{picture}(5.125,2.5)
\put(0,0){\resizebox{!}{2.5in}{
\begin{picture}(6.15,3.0) 
\put(0,0){\resizebox{3in}{!}{\includegraphics{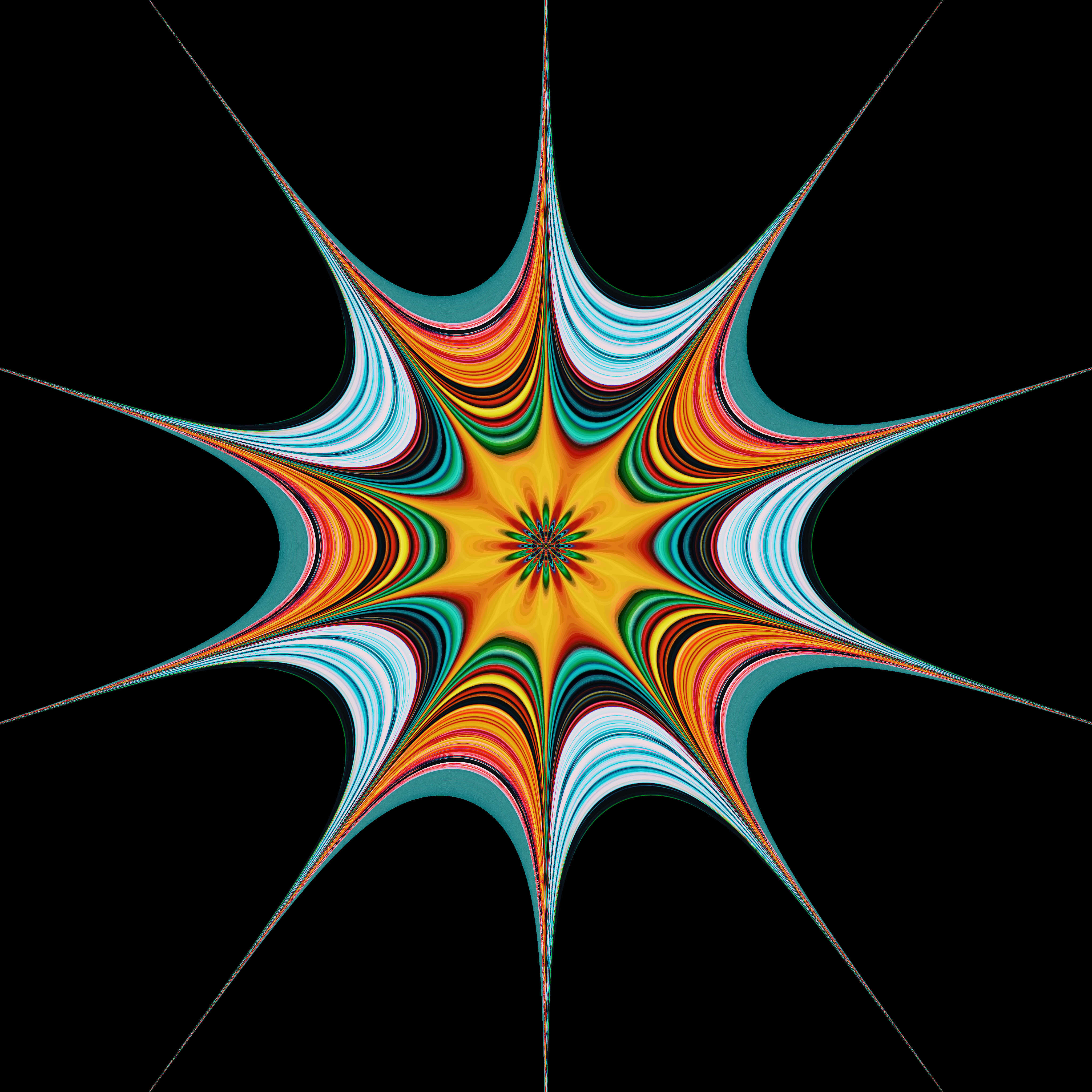}}}
\put(3.15,0){\resizebox{3in}{!}{\includegraphics{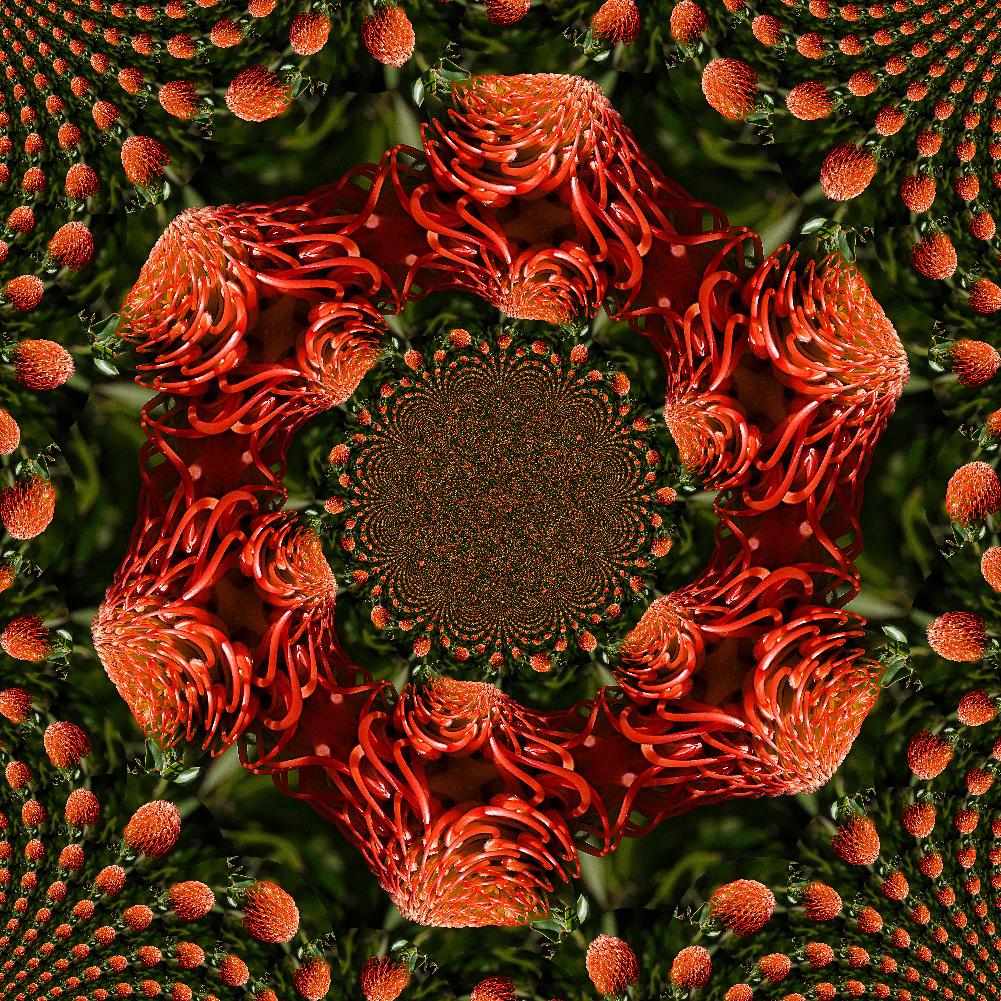}}}
\end{picture}}}
\end{picture}
\caption{\small Left: Rosette having $5$-fold rotational symmetry and reflection symmetries. Right: Rectangular view of rosette design with $6$-fold rotational symmetry plus curved lattice symmetry.\label{fig:RosetteTesselatedDesign5FoldMirroredRosette}}
\end{figure}

On the left of \Fig~\ref{fig:RosetteTesselatedDesign5FoldMirroredRosette}, we show an example of a design having
$5$-fold rotational symmetry. It also has symmetry with respect to reflection through the $x$-axis. The color map used to create the image is the \emph{Reptile} image shown in \Fig~\ref{fig:SourceImages}. The function we used to
create the design is
\begin{equation}\label{eq:5FoldDesignFunction}
f(z) = (2+3i) \left(z^5 + \overline{z}^5\right) + i\left(z^6\,\overline{z}^1 + \overline{z}^6 z^1\right)
+ \frac{i}{2000}\,\left(z^4\,\overline{z}^{-6} + z^{-6}\, \overline{z}^{-4}\right)
\end{equation}
The terms of $f(z)$ are grouped so that $f(\overline{z}) = f(z)$ clearly holds. Therefore, $f$ is symmetric with
respect to the $x$-axis reflection $R_x$. To verify that $f$ is symmetric with respect to the rotation $\rho_{2\pi/5} $, we note that $f$ is an example of a \Emph{finite} sum of the form:
\begin{equation}\label{eq:GeneralNfoldSymmetricFunction}
f(z) = \sum_{\hbox{\scriptsize$\begin{matrix}m,n\\ m\equiv n\Mod 5\end{matrix}$}} \!\! a_{m,n} z^m\, \overline{z}^n
\end{equation}
which  is symmetric with respect to $5$-fold rotation about the origin. The relation between Equations~\eqref{eq:5FoldDesignFunction} and \eqref{eq:GeneralNfoldSymmetricFunction} is that we can insure symmetry with respect to reflection through the $x$-axis by requiring $a_{m,n} = a_{n,m}$. This condition on the coefficients, $\{a_{m,n}\}$, follows by considering $f(z) = f(\overline{z})$ in terms of the coefficients of the basis functions $z^m\overline{z}^n$ in Equation~\eqref{eq:GeneralNfoldSymmetricFunction}.

\begin{figure}[!htb]\centering
\setlength{\unitlength}{1in}
\begin{picture}(5.125,2.5)
\put(0,0){\resizebox{!}{2.5in}{
\begin{picture}(6.15,3.0) 
\put(0,0){\resizebox{3in}{!}{\includegraphics{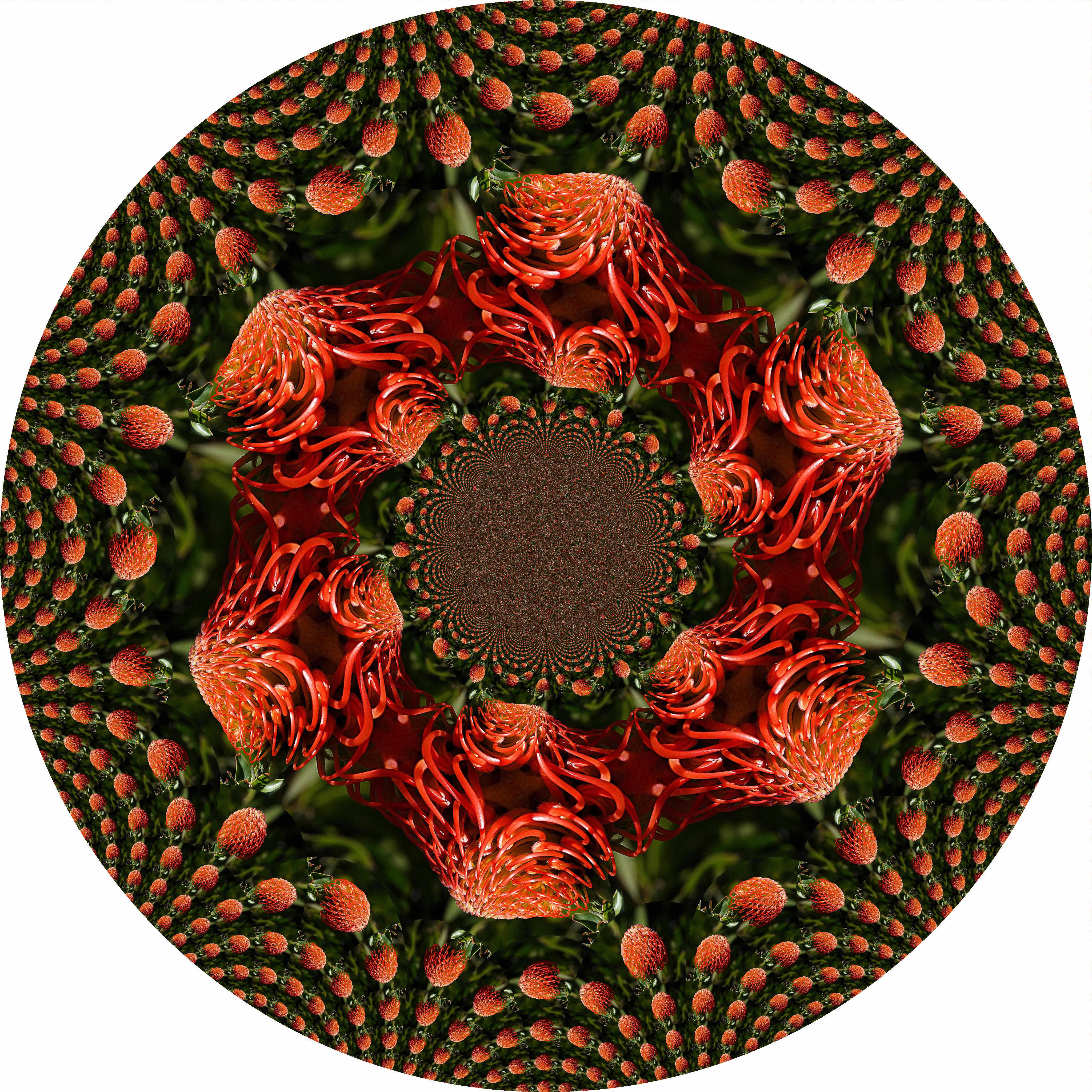}}}
\put(3.15,0){\resizebox{3in}{!}{\includegraphics{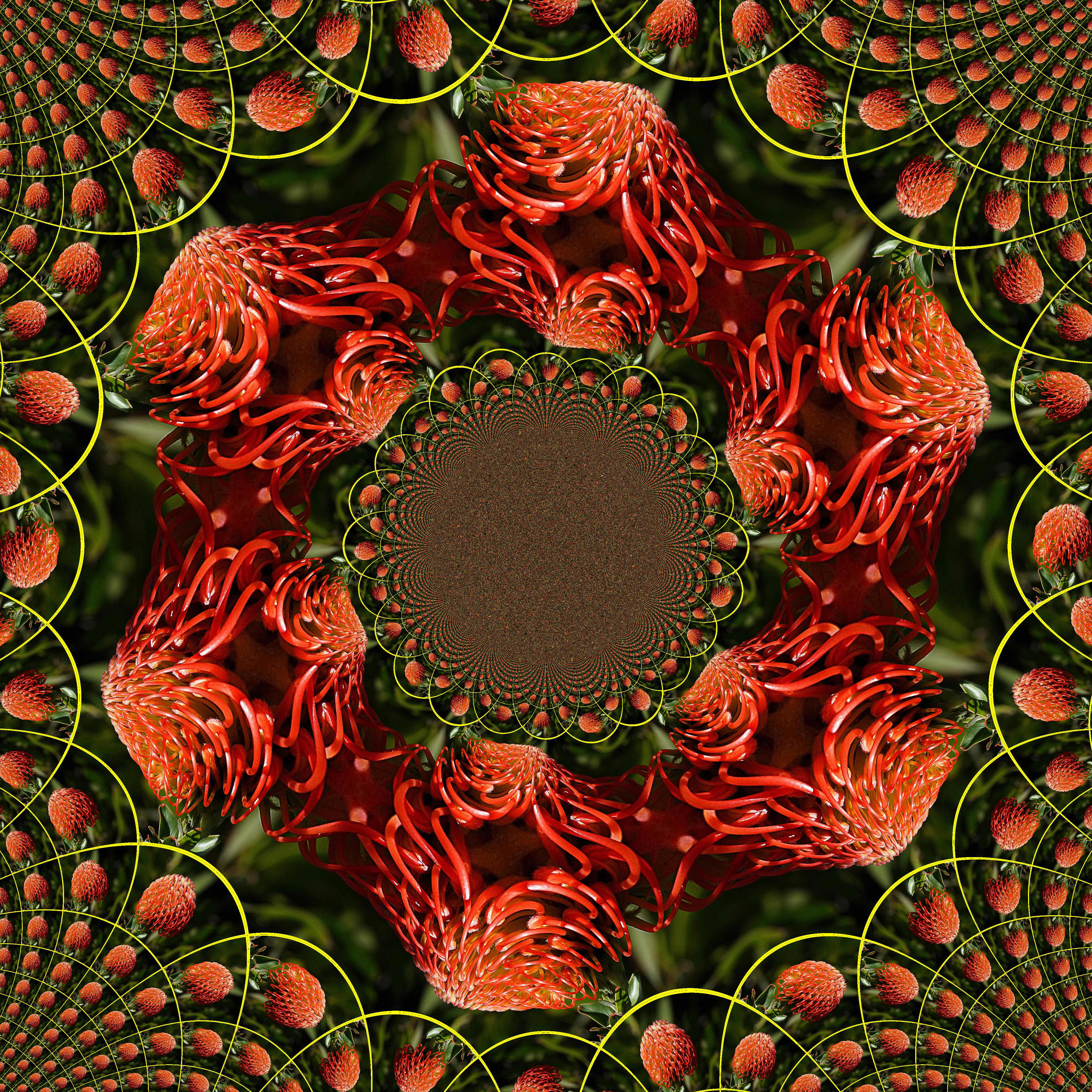}}}
\end{picture}}}
\end{picture}
\caption{\small Left: Circular view of rosette design with $6$-fold rotational symmetry plus curved lattice symmetry. Right: Curved lattice lines marked on rosette design with curved lattice symmetry. 
\label{fig:RosetteMedalionAndCurvedLatticeShownOnModRosette}}
\end{figure}

On the right of \Fig~\ref{fig:RosetteTesselatedDesign5FoldMirroredRosette}, we show another design which combines rotational symmetry with the mod operation used in Section~\ref{subsec:LatticeBasedSymmetry}. The function we used is
the following composition:
\begin{equation}\label{eq:RosetteModEquation}
z \Xrightarrow{f} (1+i) + \frac{i}{4}\, z^6 + z^{-6} = u + iv \Xrightarrow{g}
[u \Mod 1000] + i[v \Mod 1000]
\end{equation}
We are applying the function $g$ to the square image of the rosette design on the right of \Fig~\ref{fig:TesselationAndRosetteDesigns}. This is, in effect, using an infinite version of the
tessellation design in \Fig~\ref{fig:TesselationAndRosetteDesigns} as color map. (This is the
alternative way, mentioned above, for coloring pixels lying outside of a finite color map.)
The composite function satisfies $g\Bigl(f\bigl(e^{i2\pi /6} z\bigr)\Bigr) = g\Bigl(f\bigl(z\bigr)\Bigr)$.
Hence, the design on the right of \Fig~\ref{fig:RosetteTesselatedDesign5FoldMirroredRosette} has
$6$-fold rotational symmetry. This design also has several interesting features from both artistic and
mathematical perspectives:
\begin{enumerate}
\item For large magnitude $z$, we have $f(z) \sim iz^6 / 4$. Consequently, when composed with $g$, the lattice
of images of the \emph{Waratah} flower are (approximately) pre-images of $i z^6/4$ acting on the cells of a
rectangular lattice of \emph{Waratah} flower images (as shown on the left of \Fig~\ref{fig:TesselationAndRosetteDesigns}
but without the yellow grid lines). That explains the rotated, shrunken appearance of the flowers along a curved
lattice. These flowers are most easily visible at the four corners of the square image. The curved grid of the
lattice is shown clearly on the right of \Fig~\ref{fig:RosetteMedalionAndCurvedLatticeShownOnModRosette}, where
we included the yellow borders in the \emph{Waratah flower} image for our color map.
\item The square framing of the design, along with the appearance of multiple small flowers approaching the four corners
of the frame, draws the viewer's attention to the $4$-fold rotational symmetry of the design's frame. This may lead to an
ambiguity in the viewer's mind as to the validity of $6$-fold rotational symmetry for the design itself. From an
artistic perspective, we like this ambiguity. In fact, the $6$-fold symmetry only holds for points in the image that
remain within the image boundaries upon rotation by $2\pi/6$. On the left of \Fig~\ref{fig:RosetteMedalionAndCurvedLatticeShownOnModRosette} we show a larger scale view of the design,
within a circular frame. This latter design clearly displays $6$-fold rotational symmetry.
\item For small magnitude $z$, we have $f(z) \sim z^{-6}$. Consequently, the curved lattice contains an
infinity of pre-images approaching the pole of multiplicity $6$ at $z = 0$. The sizes of the pre-images decrease
rapidly to $0$ as $z\to 0$, hence they become incapable of realization with the digital images that we are working
with here. This explains the random looking colored pixels lying within a region interior to the six flower-like
parts of image. Those flower-like parts lie across the unit-circle, where $|z| \approx 1$. Some of the curved lattice
within the unit-circle is visible in the image on the right of \Fig~\ref{fig:RosetteMedalionAndCurvedLatticeShownOnModRosette}.
\end{enumerate}
\noindent
\subsubsection{Animations\label{subsubsec:Animations}}
We have also created animations involving our symmetric designs. We will just give one illustration of the ideas involved. The animation we describe can be found at this link:
\begin{verbatim}
https://www.youtube.com/watch?v=Pf1vJPywXWs
\end{verbatim}
We produced this animation in the following way. First, we selected a unit-length complex number $p = e^{i\theta}$ for small
positive $\theta$.  We then created a succession of designs using
the mappings:
\begin{equation}\label{eq:6FoldSymmetryDesignAnimation}
z \Xrightarrow{} (1 + i) + p^n\, (i/4)z^6 + \frac{1}{p^n}\, \frac{1}{z^6} \; = u + i v\;
\Xrightarrow{\mathcal{C}}\; [u\Mod 1000] + i [v\Mod 1000]
\end{equation}
for $n = 0, 1, 2, \dots, 99$. This produced $100$ separate designs that are displayed one after another to create
the animation. The interesting thing about this animation is that it will appear to rotate in the positive angular
direction on its outer part versus a negative angular rotation on its inner part. The reason for this is that $p^n$ has angle $n\theta$ for $n = 0, 1, 2, \dots, 99$, which is a sequence of increasing
angles in the counter-clockwise direction. The outer part of the design corresponds to values
of $z$ that have lengths larger than $1$, and for larger values of $|z|$ we have
\[
(1 + i) + p^n\, (i/4)z^6 + \frac{1}{p^n}\, \frac{1}{z^6} \approx p^n\, (i/4)z^6
\]
so the outer part of the design successively rotates through the angles $n\theta = \arg(p^n)$, $n=0, 1, 2, \dots, 99$. A similar argument shows that for $z$ close to $0$, the design simultaneously rotates through the angles $-n\theta$.
\subsection{Rotationally Symmetric Wallpaper Designs\label{subsec:WallpaperSymmetry}}
\subsubsection{Four-fold and Two-fold Wallpaper Symmetry\label{subsubsec:FourFoldAndTwoFoldWallpaperSymmetry}}
\begin{figure}[!htb]\centering
\setlength{\unitlength}{1in}
\begin{picture}(5.852,2.35)
\put(0,0){\resizebox{!}{2.2in}{
\begin{picture}(6.65,2.5)
\put(0,0){\includegraphics[width=2.7in]{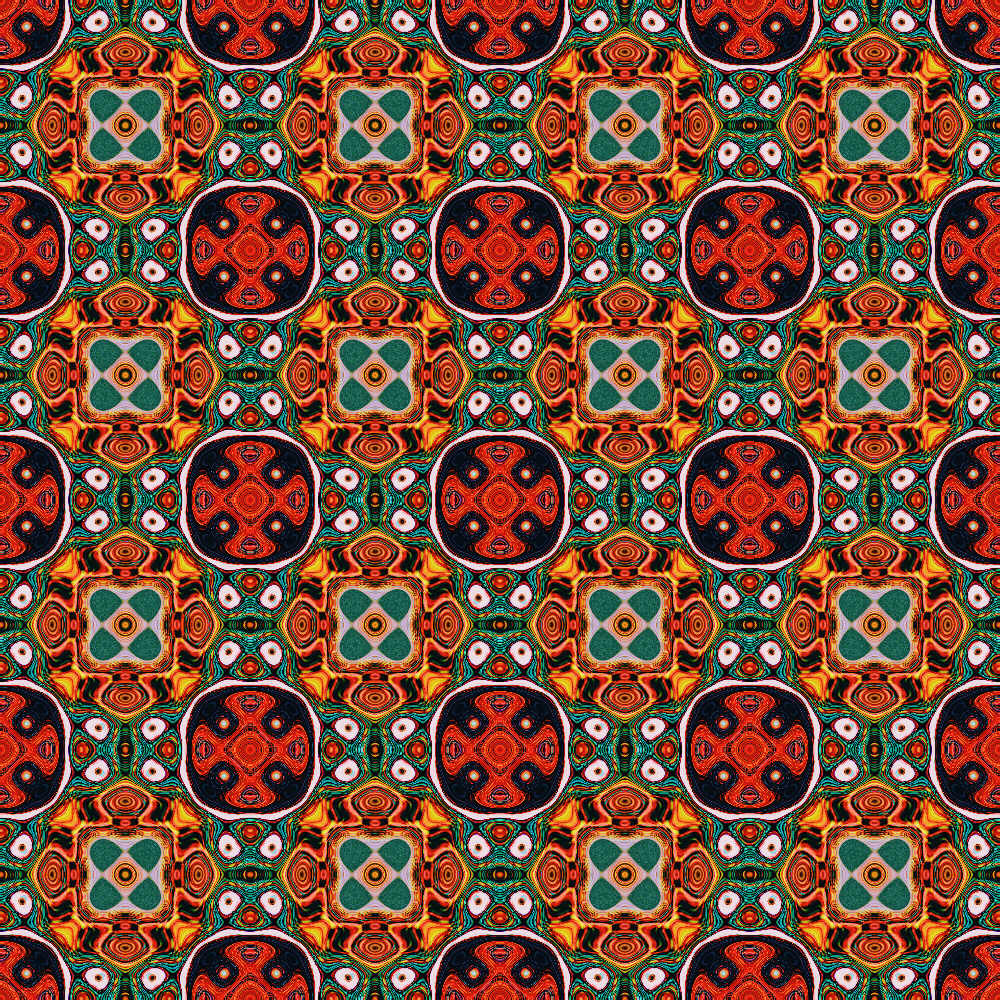}}
\put(2.8,0.95){\includegraphics[width=0.8262in]{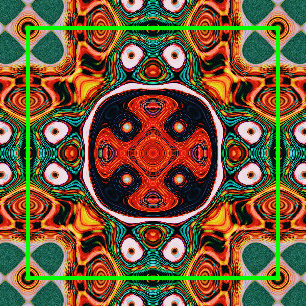}}
\put(3.95,0){\includegraphics[width=2.7in]{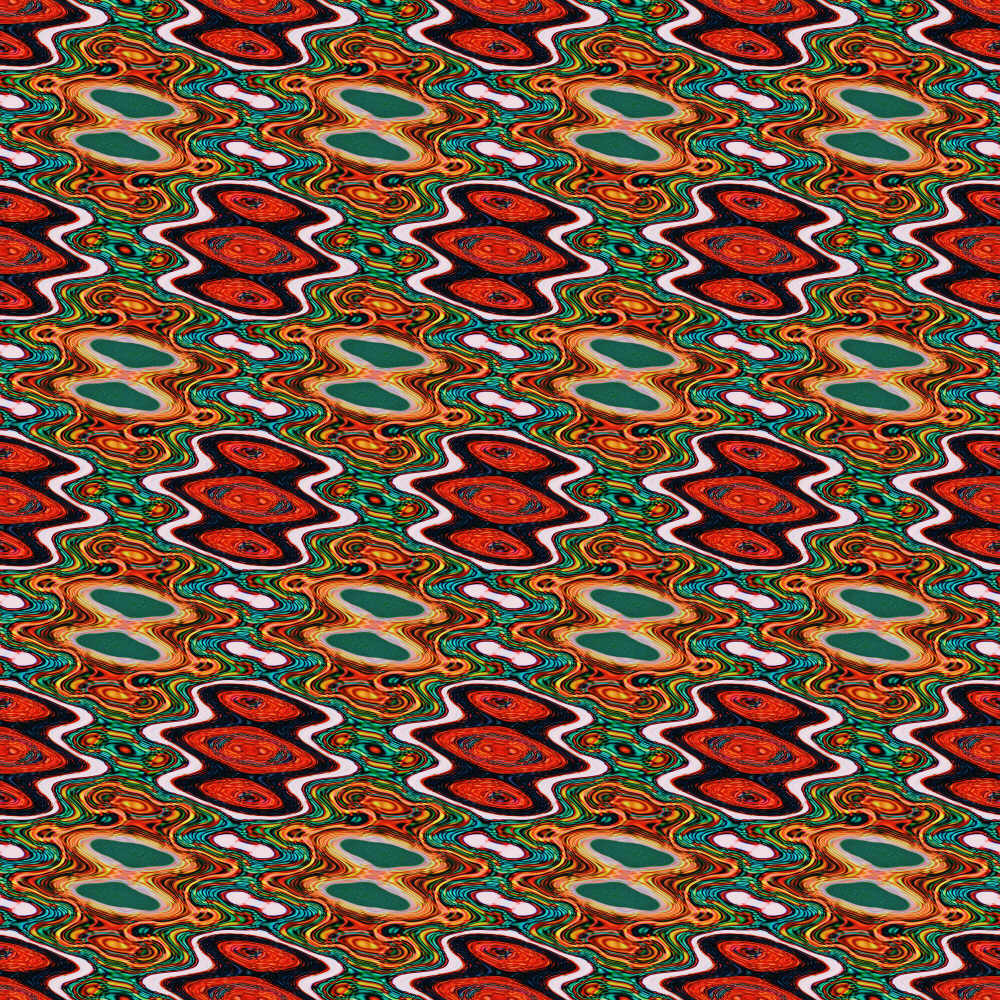}}
\end{picture}}}
\end{picture}
\caption{\small Left: $4$-fold rotationally symmetric design, with additional symmetries.
Unit cell on its right marked in green. Right: $2$-fold rotationally symmetric design with no additional symmetries.\label{fig:4foldSymmetricWallpaper}}
\end{figure}

In \Fig~\ref{fig:4foldSymmetricWallpaper} we show a wallpaper design that is symmetric with respect to
$4$-fold symmetry and various reflection symmetries. For this design, we used the \emph{Reptile} image.
To create a function $f$ having square lattice symmetry and four-fold rotational symmetry, we follow the
procedure described by Farris~\cite{Farris2015}. First, we use as a basis the set of complex exponentials $\{E_{m,n}(z) = e^{2\pi i(mx+ny)}\}$ for all $m, n \in\mathbb{Z}$ and $z = x + iy \in \mathbb{C}$. Any finite, or convergent, sum $\sum a_{m,n}E_{m,n}(z)$ is guaranteed to have the required translational symmetry. To obtain rotational symmetry, we employ group averaging. In this case, the group is the rotations given by powers of $i = e^{2\pi i/4}$. The group average $W_{m,n}$ of $E_{m,n}$ is defined by  $W_{m,n}(z) = \frac{1}{4}\sum^3_{k=0} E_{m,n}(i^k z)$.
Any finite, or convergent, linear combination $f(z) = \sum a_{m,n} W_{m,n}(z)$ is guaranteed to have both square lattice symmetry and four-fold rotational symmetry about the origin.

The specific function $f$ that we employed is
\begin{equation}\label{eq:4FoldSymmetryFunction}
\begin{split}
f(z) &= \left[W_{1,0}(z) + W_{0,-1}(z)\right] + 0.5\left[W_{1,5}(z) + W_{-5,-1}(z)\right] \\
&\qquad + 0.1i\left[W_{-2,4}(z) + W_{-4,2}(z)\right] - 0.05i\left[W_{-6,3}(z) + W_{-3,6}(z)\right]
\end{split}
\end{equation}
In this equation, the functions $\{W_{m,n}(z)\}$ are grouped in order to enforce reflectional symmetry through the $x$-axis. To see this, we first rewrite the terms of $W_{m,n}$ to remove the powers of $i$. The term $E_{m,n}(iz)$ satisfies
\begin{align*}
E_{m,n}(iz) &= E_{m,n}(-y + ix)\\
&= e^{i2\pi (-my + nx)}\\
&= E_{n,-m}(z)
\end{align*}
Iterating this relation yields $E_{m,n}(i^2 z) = E_{-m,-n}(z)$ and
$E_{m,n}(i^3 z) = E_{-n,m}(z)$. Therefore, 
\begin{equation}\label{eq:4FoldGroupAverageFormula}
W_{m,n}(z) = \frac{1}{4}\,\left[E_{m,n}(z) + E_{n,-m}(z) + E_{-m,-n}(z) +  E_{-n,m}(z)\right]
\end{equation}
Now, if we apply reflection through the $x$-axis to $E_{m,n}(z)$, we obtain
\begin{align*}
E_{m,n}(\overline{z}) &= E_{m,n}(x - iy)\\
&= E_{m,-n}(z)
\end{align*}
Iterating this relation, we obtain
$E_{n,-m}(\overline{z}) = E_{n,m}(z)$, $E_{-m,-n}(\overline{z}) = E_{-m,n}(z)$, and
$E_{-n,m}(\overline{z}) = E_{-n,-m}(z)$. Therefore,
\begin{align*}
W_{m,n}(\overline{z}) &= \frac{1}{4}\,\left[E_{m,-n}(z) + E_{n,m}(z) + E_{-m,n}(z) +  E_{-n,-m}(z)\right]\\
&= W_{-n,-m}(z)
\end{align*}
Consequently, terms of the form
\[
a\left[W_{m,n}(z) + W_{-n,-m}(z)\right]
\]
appearing in Equation~\eqref{eq:4FoldSymmetryFunction}, are symmetric with respect to reflection through the $x$-axis.

The symmetry with respect to $R_x$ is apparent
in the design in \Fig~\ref{fig:4foldSymmetricWallpaper}, but there are other symmetries as well. For example, it has
reflectional symmetry $R_y$ about the $y$-axis. Since the $y$-axis is a rotation by $2\pi/4$ of the $x$-axis, the function $f$ must be symmetric with respect to $R_y$. More precisely, $R_y$ equals the conjugation operation
$\rho_{2\pi/4}^{-1} R_x\, \rho_{2\pi/4}$ in the symmetry group $\S_f$. Other reflection symmetries follow from
group operations, such as reflection through the line $y = x$ and reflection through the line $y = -x$.

This wallpaper design is symmetric with respect to translation by the two independent vectors $1$ and $i$ in $\C$. The cells of the lattice generated by $1$ and $i$ are clearly evident. To the right of the design, we show a \Emph{unit cell}, the contents of which generate the entire design via translations by $m\cdot 1 + n\cdot i$, for $m, n \in \Z$.  One interesting feature of the design is that the center of the unit cell is a point of $4$-fold rotational symmetry for the design. This can be proved as follows. Assuming this unit cell has its lower left corner at the origin, its center is $\frac{1}{2} + \frac{1}{2}\, i$. The following mappings
\[
\frac{1}{2} + \frac{1}{2}\, i \Xrightarrow{\cdot i} \frac{-1}{2} + \frac{1}{2}\, i \Xrightarrow{\tau_1} \frac{1}{2} + \frac{1}{2}\, i
\]
imply that $\tau_1 \circ \rho_{2\pi /4}$ preserves $\frac{1}{2} + \frac{1}{2}\, i$ and rotates the unit cell by $2\pi / 4$. Consequently, $\frac{1}{2} + \frac{1}{2}\, i$ is a center of $4$-fold rotational symmetry.
A similar argument shows that the midpoint of the top side of the unit cell is a center of $2$-fold rotational symmetry.
Hence, by four-fold rotation, the midpoints of each side of the unit cell are centers of $2$-fold symmetry.

Two-fold rotationally symmetric wallpaper patterns can be generated in multiple ways. One, rather elementary, way is to
just stretch a $4$-fold pattern to create a rectangular lattice.
The rectangular lattice for a $2$-fold design, that is not $4$-fold, is generated by
the basis vectors $1$ and $ri$, where $r > 0$ and $r\neq 1$. Such a design can be generated from a $4$-fold design
by simply stretching or shrinking the $4$-fold design in the vertical direction. With the digital images we are
creating, that can be done by even the most rudimentary image processing programs. A second method would be to
use a group average approach, as we did with the $4$-fold case. For instance, we can use group
averages $W_{m,n}(z) = \frac{1}{2}\bigl(E_{m,n}(z) + E_{m,n}(-z)\bigr)$. A finite, or convergent, linear combination
$f(z) = \sum a_{m,n} W_{m,n}(z)$ will then have both $2$-fold rotational and square lattice symmetry. An example
is shown on the right of \Fig~\ref{fig:4foldSymmetricWallpaper}. The function $f(z)$ that we used for this design
is
\begin{equation}\label{eq:2FoldWallpaperEquation}
f(z) = W_{1,0}(z) + W_{0,-1}(z) + 0.5 W_{1,5}(z) + 0.1 i W_{-2,4}(z) - 0.05 i W_{-6,3}(z)
\end{equation}
Unlike the design on the left of \Fig~\ref{fig:4foldSymmetricWallpaper}, this design has no additional symmetries.
This corresponds to the absence of grouping of related terms in Equation~\eqref{eq:2FoldWallpaperEquation}, of the kind
that we have in Equation~\eqref{eq:4FoldSymmetryFunction}.
\subsubsection{Three-Fold Wallpaper Symmetry\label{subsubsec:3FoldWallpaperSymmetry}}
\begin{figure}[!htb]\centering
\setlength{\unitlength}{1in}
\begin{picture}(6.5,2.55)
\put(0,1.8){\includegraphics[width=1.2in]{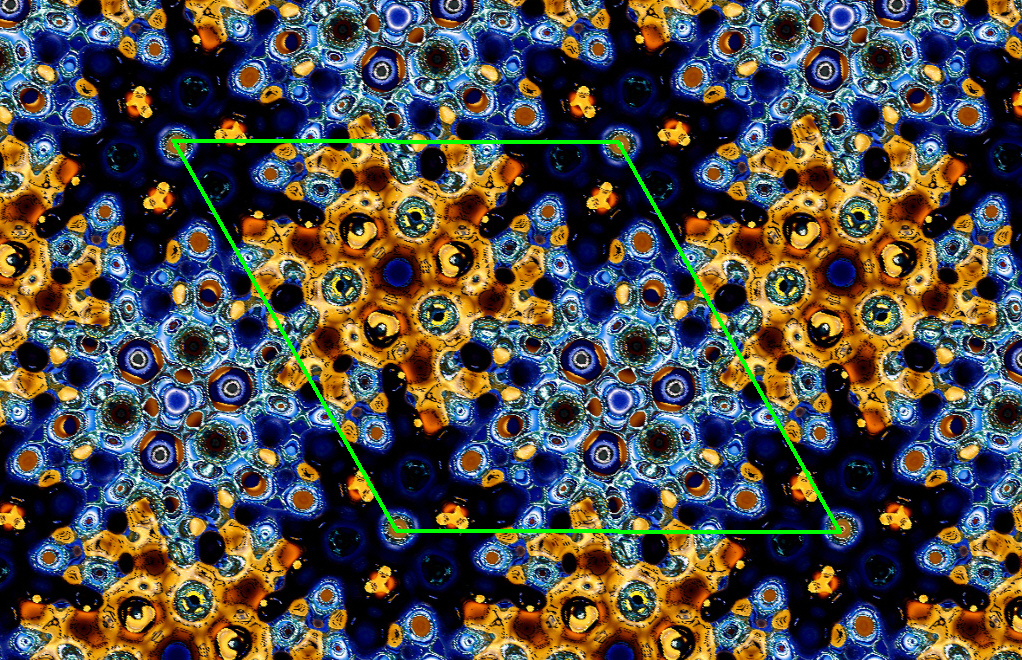}}
\put(0,0){\includegraphics[width=2.4in]{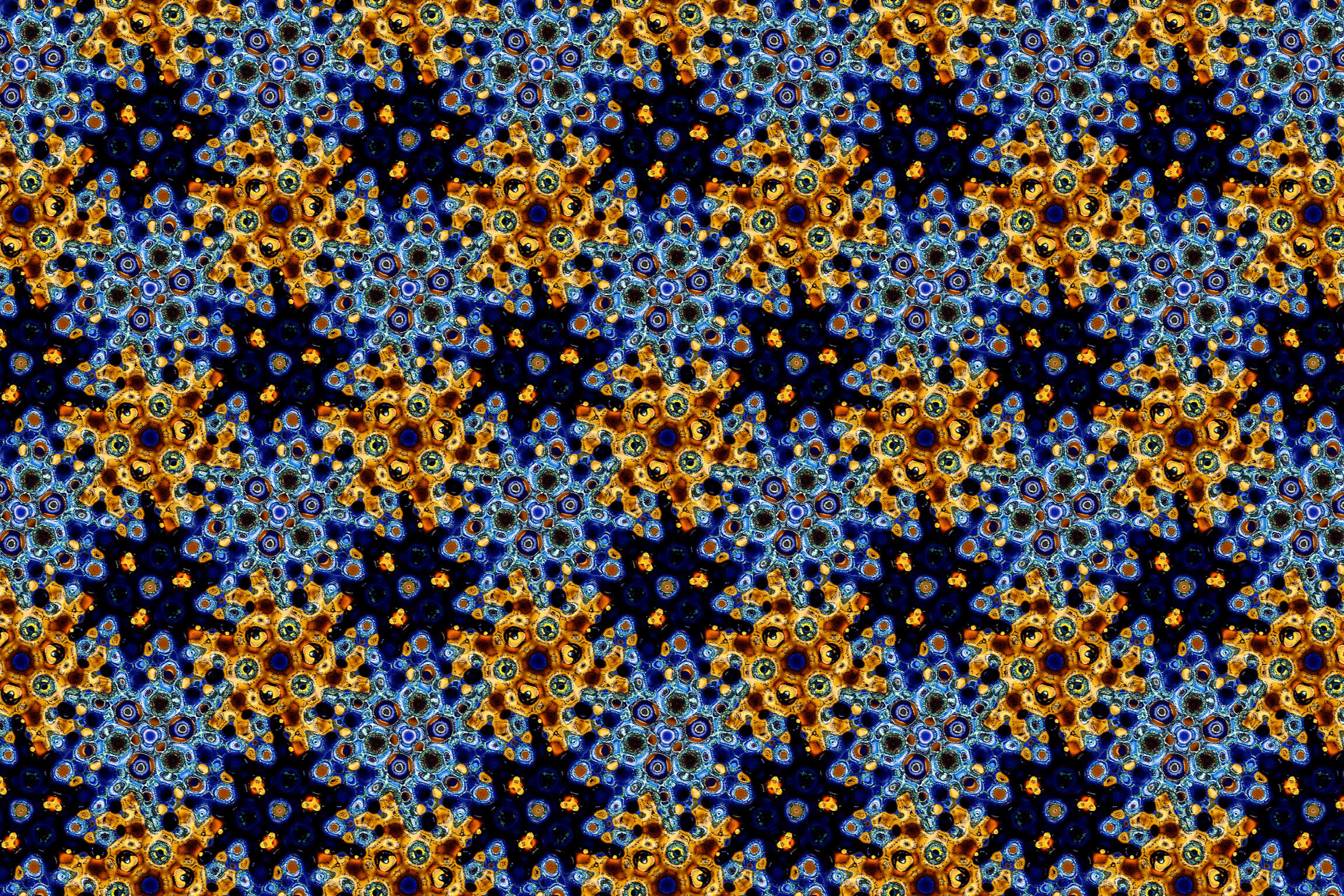}}
\put(2.48,0){\includegraphics[width=4in]{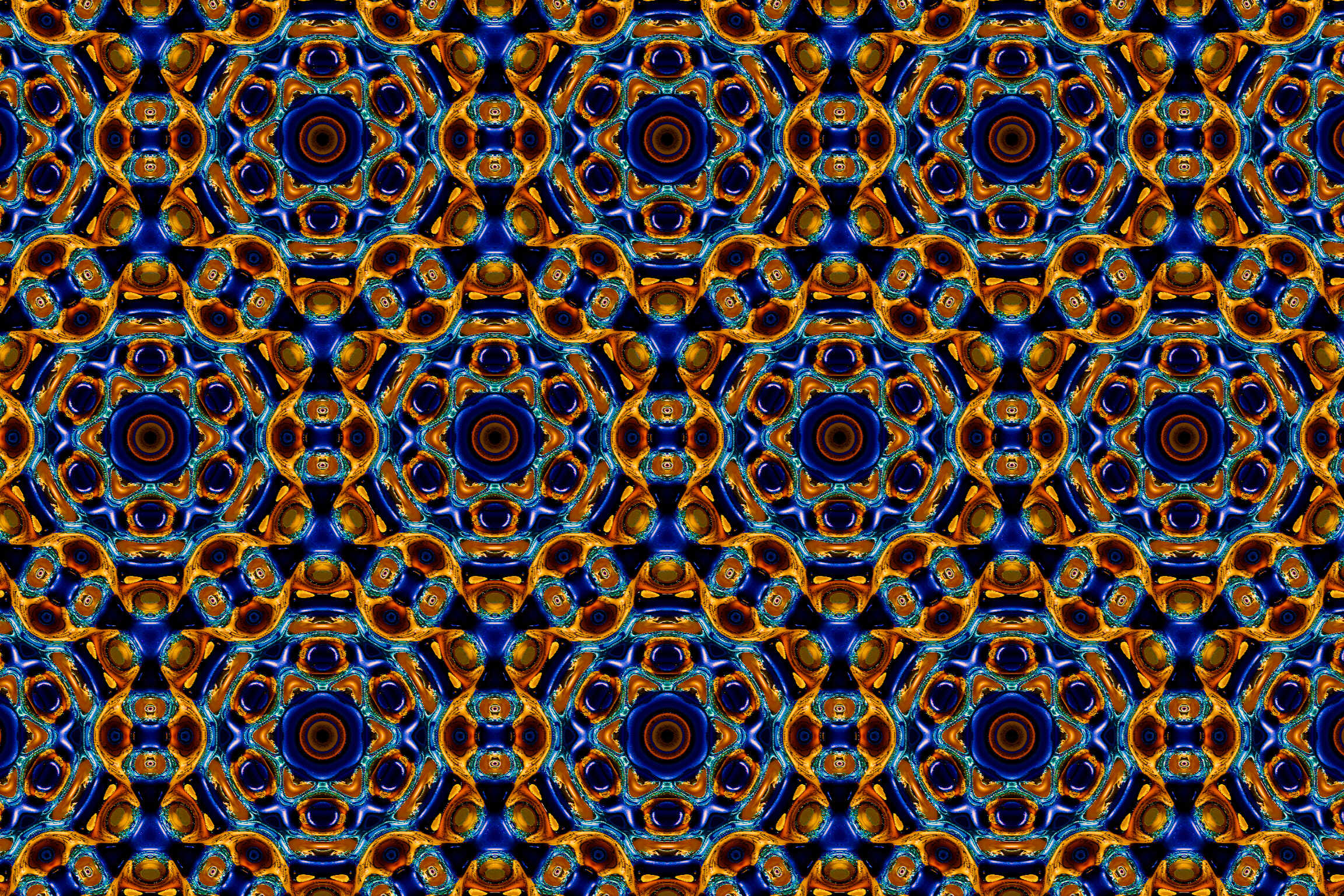}}
\end{picture}
\caption{\small Bottom Left: Wallpaper design with $3$-fold rotational symmetry. Top Left: a rhombic unit cell, marked in green.
Right: Wallpaper design with $6$-fold rotational symmetry and additional reflection symmetries. \label{fig:3and6foldSymmetricWallpaper}}
\end{figure}

On the left of \Fig~\ref{fig:3and6foldSymmetricWallpaper} we show a wallpaper design that is symmetric with respect to
$3$-fold symmetry.
Above the design is an image of a unit cell for the lattice of the design, a rhombus with sides constructed from the complex numbers $1$ and $\omega = e^{i2\pi/3}$. The rhombic lattice for translational symmetry is defined by the vectors $m + n\omega$, for all $m, n\in\Z$ and $\omega = e^{i2\pi/3}$. 
Every complex number $z$ can be expanded uniquely as $z = u + v\omega$ for unique $u, v\in\R$. The functions, $E_{m,n}(z)$ for $m, n\in\Z$, are defined by $E_{m,n}(z) = E_{m,n}(u + v\omega) = e^{i2\pi (mu+nv)}$ for all $u, v\in \R$. These functions are periodic over the cells of the rhombic lattice generated by $1$ and $\omega$. The group averages, $\{W_{m,n}\}$, for the group of $3$-fold rotations about the origin, are defined by
\begin{equation}
W_{m,n}(z) = \frac{1}{3}\left[E_{m,n}(z) + E_{m,n}(\omega z) + E_{m,n}(\omega^2 z)\right]
\end{equation}
This definition of $W_{m,n}(z)$ ensures that it has $3$-fold rotational symmetry. To verify that it has translational symmetry,
we need to further examine the terms $E_{m,n}(\omega z)$ and $E_{m,n}(\omega^2 z)$. Since
$z = \omega$ satisfies the factored equation $0 = (z^3 - 1) = (z-1)(z^2 + z + 1)$, it follows that
$\omega^2 = -1 - \omega$. Therefore,
\begin{align*}
E_{m,n}(\omega z) &= E_{m,n}(u\omega + v\omega^2)\\
&= E_{m,n}\bigl(-v + (u-v)\omega\bigr)\\
&= e^{i2\pi \bigl(-mv + n(u-v)\bigr)}\\
&= e^{i2\pi \bigl(nu + (-m-n)v\bigr)}\\
&= E_{n, -m-n}(z)
\end{align*}
Thus, $E_{m,n}(\omega z) = E_{n, -m-n}(z)$. Iterating this relation,
we obtain $E_{m,n}(\omega^2 z) = E_{-m-n,m}(z)$. These relations show that
\begin{equation}\label{eq:BasisExpansionOfGroupAverageFunction}
W_{m,n}(z) = \frac{1}{3}\left[E_{m,n}(z) + E_{n,-m-n}(z) + E_{-m-n,m}(z)\right]
\end{equation}
Hence, $W_{m,n}$ is a linear combination of the basis functions $\{E_{m,n}\}$, so it also enjoys
translational symmetry over the rhombic lattice. Consequently any function $f$ defined by a
finite, or convergent, sum $\sum_{m,n} a_{m,n} W_{m,n}(z)$ has both $3$-fold rotational symmetry and
rhombic lattice symmetry. Furthermore, because of the group operation $E_{m,n}\cdot E_{j,k} = E_{m+j,n+k}$
for the basis functions $\{E_{m,n}\}$, we can also include products of the form $W_{m,n}\cdot W_{j,k}$ in
the terms for $f$. Specifically, the function we used to create the design on the bottom left of \Fig~\ref{fig:3and6foldSymmetricWallpaper} is
\[
f(z) = 2W_{2,3}(z)\cdot W_{1,4}(z) + W_{1,0}(z)
\]
It has $3$-fold rotational symmetry because each of the factors, $W_{2,3}$ and $W_{1,4}$, and the term $W_{1,0}$ have that symmetry. Moreover, it has translational symmetry over the lattice because $f$ is a linear combination of the
basis functions $\{E_{m,n}\}$.
\begin{small}\begin{remark}\label{rem:FortunateCaseOf3FoldSymmetry}
In this example, we were indeed fortunate that a rotation by $2\pi/3$ maps each basis function
$E_{m,n}$, periodic over the lattice, into another such basis function. That fact allowed us to create a $3$-fold rotationally symmetric wallpaper design. The \textbf{Crystallographic Restriction Theorem} below implies that this can not happen for most $n$-fold rotational symmetries.
\end{remark}\end{small}
\subsubsection{Six-Fold Wallpaper Symmetry\label{subsubsec:6FoldWallpaperSymmetry}}
In \Fig~\ref{fig:3and6foldSymmetricWallpaper} we show a wallpaper design that is symmetric with respect to
$6$-fold symmetry and various reflection symmetries. 
Before we discuss the mathematics of creating this design, we will take a moment to comment on its artistic features.
The design presents an appearance of interlocking circular regions. The larger regions have a central, three-armed
cross that appears to be a center of three-fold rotational symmetry. Slightly smaller circular regions have
dark blue, hexagonally shaped figures about their centers, which appear to have $6$-fold rotational symmetry. If
you relax your focus slightly, the whole design appears to float in the background as your attention rapidly shifts
between these interlocking circles with different symmetries. These features slightly disguise the location of the
rhombus shaped cells that form the lattice for the design.

We used the \emph{Buoy} image for the color map of this design. To create a function $f$ having rhombic lattice symmetry and six-fold rotational symmetry, we note that $6$-fold rotational symmetry about the origin is equivalent to $3$-fold rotational symmetry combined with $2$-fold rotational symmetry. Two-fold rotational symmetry about the origin in $\C$
corresponds to the mapping $\rho_\pi: z \to -z$. Since $E_{m,n}(-z) = E_{-m,-n}(z)$, it follows that
$W_{m,n}(-z) = W_{-m,-n}(z)$ for all $m,n\in\Z$. Consequently, if $f$ is a finite, or convergent, sum of terms
of the form
\[
a\bigl[W_{m,n}(z) + W_{-m,-n}(z)\bigr]
\]
then $f$ is symmetric with respect to $R_\pi$ and $R_{2\pi/3}$ and thus symmetric with respect to $R_{2\pi/6}$.

The function $f$ that we used to create the design in \Fig~\ref{fig:3and6foldSymmetricWallpaper} is somewhat complicated.
Suffice it to say that it had the following form:
\begin{equation}\label{eq:6FoldSymmetryFunction}
\begin{split}
f(z) &= a\left[W_{2,3}(z) + W_{3,2}(z)\right] + a\left[W_{-2,-3}(z) + W_{-3,-2}(z)\right] \\
&\qquad + b\left[W_{1,5}(z) + W_{5,1}(z)\right] + b\left[W_{-1,-5}(z) + W_{-5,-1}(z)\right]\\
&\qquad + c\left[W_{3,4}(z) + W_{4,3}(z)\right] + c\left[W_{-3,-4}(z) + W_{-4,-3}(z)\right]\\
\end{split}
\end{equation}
for certain complex constants, $a$, $b$, and $c$. The terms for $f(z)$ can be regrouped so that $f$ is
a linear combination of terms of the form $W_{m,n}(z) + W_{-m,-n}(z)$. Hence, $f$ has $6$-fold rotational symmetry.
However, the terms for $f$ were grouped in the pairs shown in Equation~\eqref{eq:6FoldSymmetryFunction} in order to enforce
reflection symmetry $R_x$ about the $x$-axis. To see that this symmetry holds, we examine the effect of
$R_x: z \to \overline{z}$ on $E_{m,n}$:
\begin{align*}
E_{m,n}(\overline{z}) &= E_{m,n}(u + v\overline{\omega})\\
&= E_{m,n}\bigl(u + v(-1 - \omega)\bigr)\qquad\text{[since $\overline{\omega} = \omega^2$]}\\
&= e^{i2\pi \bigl(mu + (-m-n)v\bigr)}\\
&= E_{m,-m-n}(z)
\end{align*}
Thus, $E_{m,n}(\overline{z}) = E_{m,-m-n}(z)$. Iterating this relation, we obtain
\[
E_{n,-m-n}(\overline{z}) = E_{n,m}(z)\qquad \text{and}\qquad
E_{-m-n,m}(\overline{z}) = E_{-m-n,n}(z)
\]
Combining these relations with Equation~\eqref{eq:BasisExpansionOfGroupAverageFunction}, we obtain $W_{m,n}(\overline{z}) = W_{n,m}(z)$. Consequently, terms of the form
\[
a\left[W_{m,n}(z) + W_{n,m}(z)\right]
\]
are symmetric with respect to reflection through the $x$-axis. Since $f$ is a finite sum of such terms, it has
this symmetry, too. Our design has other symmetries as well. For example, since $\rho_\pi \circ R_x: z \to -\overline{z}$,
the function $f$ is mirror symmetric through the $y$-axis.

At the beginning of our discussion of this $6$-fold symmetric design, we mentioned circles having centers of $3$-fold
rotational symmetry within the rhombic cells of the design's lattice. The point $\frac{2}{3} + \frac{1}{3}\,\omega$ within
the rhombic unit cell is a center of $3$-fold symmetry. The following mappings
\[
\frac{2}{3} + \frac{1}{3}\,\omega \Xrightarrow{\cdot \omega} \frac{-1}{3} + \frac{1}{3}\, \omega
\Xrightarrow{\tau_1} \frac{2}{3} + \frac{1}{3}\,\omega
\]
imply that $\tau_1\circ \rho_{2\pi/3}$ preserves $\frac{2}{3} + \frac{1}{3}\,\omega$ and rotates the rhombic unit cell
by $2\pi/3$. Consequently, $\frac{2}{3} + \frac{1}{3}\,\omega$ is a center of $3$-fold rotational symmetry. Furthermore,
\[
\frac{1}{3} + \frac{2}{3}\,\omega \Xrightarrow{\cdot \omega} \frac{-2}{3} - \frac{1}{3}\, \omega
\Xrightarrow{\tau_{1+\omega}} \frac{1}{3} + \frac{2}{3}\,\omega
\]
implies that $\tau_{1+\omega}\circ \rho_{2\pi/3}$ preserves $\frac{1}{3} + \frac{2}{3}\,\omega$, and so it
is a second point of $3$-fold symmetry within the rhombic unit cell.
\subsubsection{The Crystallographic Restriction\label{subsubsec:CrystallographicRestriction}}
We have shown that wallpaper designs can be generated with $n$-fold rotational symmetry when $n=2, 3, 4$, and $6$. In fact, these are the only possible $n$-fold rotationally symmetric wallpaper designs.
\begin{theorem}[The Crystallographic Restriction] \textsl{An $n$-fold rotationally symmetric wallpaper design is only possible if
$n = 2$, $3$, $4$, or $6$.}
\end{theorem}
\begin{small}\begin{proof}
In the standard basis for $\R^2 \equiv \C$, the $n$-fold rotation $\rho_{2\pi /n}$ has matrix form
\[
\begin{pmatrix}\cos(2\pi/n) & -\sin(2\pi/n)\\ \sin(2\pi/n) & \phantom{-}\cos(2\pi/n)\end{pmatrix}
\] The trace of this matrix is $2\cos(2\pi/n)$. However, the trace is invariant under change of basis. For the basis that generates the lattice of cells for the wallpaper design, the trace of the matrix for $\rho_{2\pi/n}$ must be an
integer. Hence $2\cos(2\pi/n) = k$ for some $k\in\Z$. Therefore, we have
\begin{equation}\label{eq:CrystallographicRestriction}
\cos(2\pi/n) = k/2, \quad k\in\Z
\end{equation}
The only possible integers $k$ for which Equation~\eqref{eq:CrystallographicRestriction} can hold are
$k = 0, \pm 1, \pm 2$, which yield $n = 2, 3, 4$, and $6$.
\end{proof}\end{small}\noindent
There is an interesting analysis of this restriction when $n = 5$ in \cite{Farris2013ForbSym}. Some
Matlab$^\circledR$ programs for creating designs when $n = 5$ are in \cite[p.~98]{Terras}. These designs are related to \emph{quasicrystals}~\cite{Terras, ref:Senechal}.
\subsection{Summary\label{subsec:Summary}}
We have shown a number of symmetric designs generated by the application of complex analysis to the geometry
of the Euclidean plane. The mathematics we have used is widely employed in
crystallography~\cite{ref:XrayCrystallography, ref:Senechal}.
Many more designs, and a more thorough treatment including the relation to crystallography, can be found in the
book by Farris~\cite{Farris2015}. In the next section, we describe symmetric designs that use properties
of complex analysis in hyperbolic geometry.
\section{Designs in the Hyperbolic Upper Half-Plane\label{sec:DesignsUsingHyperbolicGeometry}}
\begin{figure}[!htb]\centering
\setlength{\unitlength}{1in}
\begin{picture}(5.33807829,5.05)
\put(0,0){\resizebox{!}{5in}{
\begin{picture}(6,5.62)
\put(0,3.66){\includegraphics[width=6in]{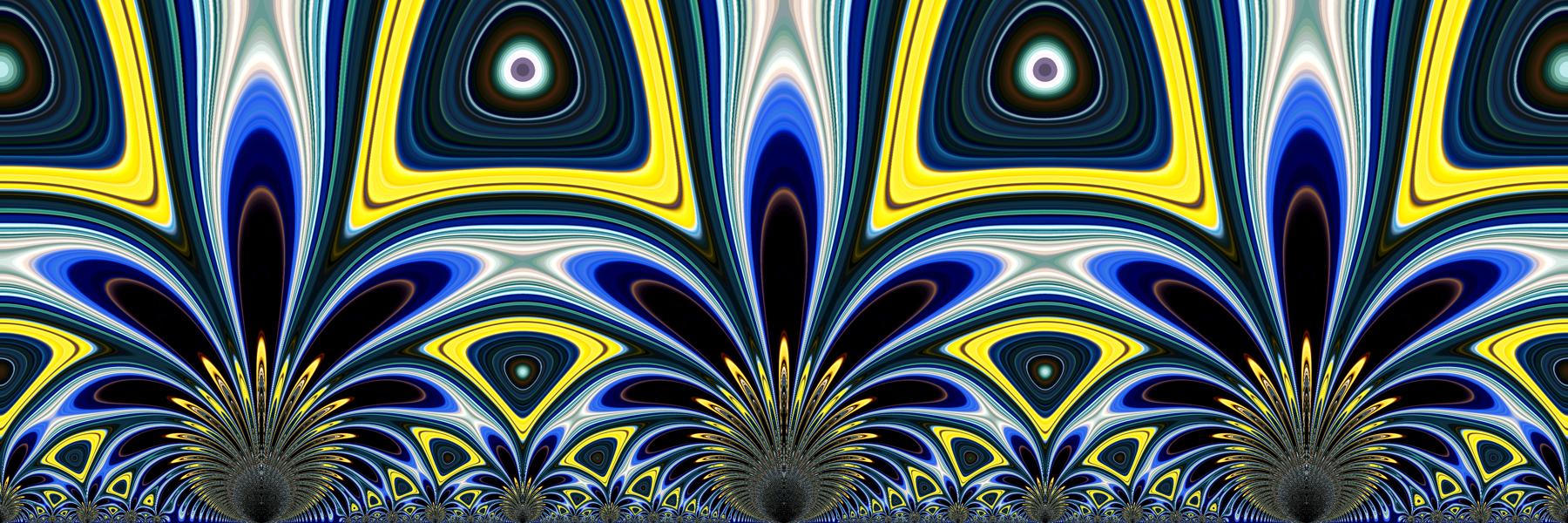}}
\put(0,1.58){\includegraphics[width=6in]{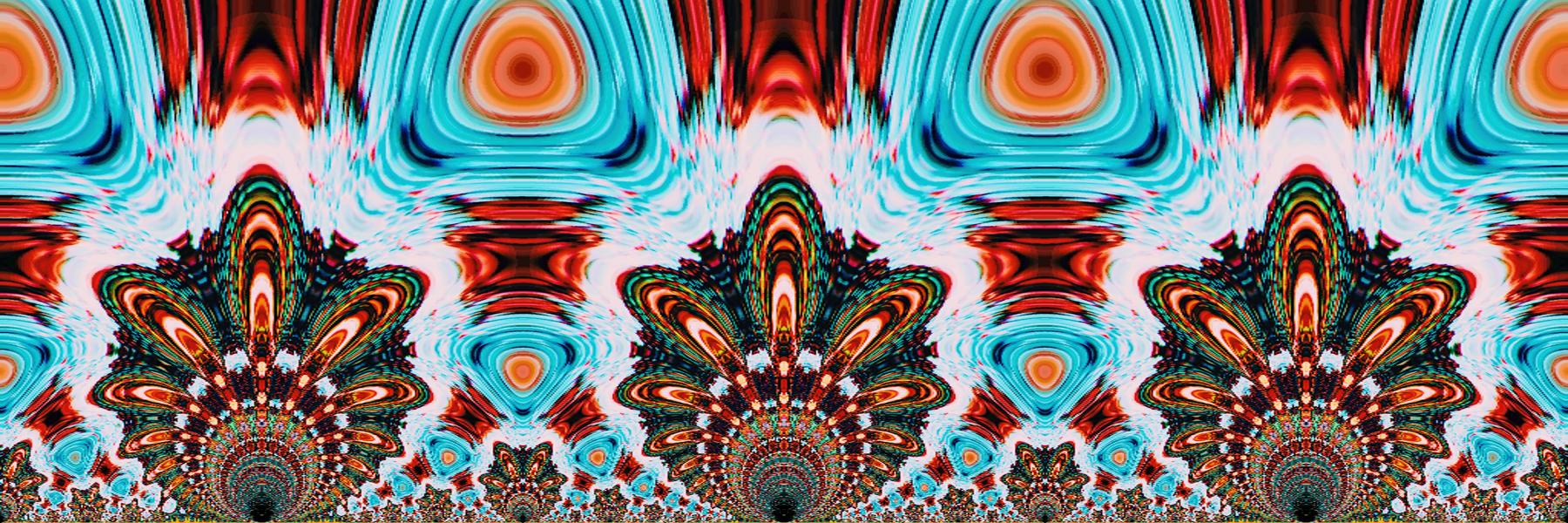}}
\put(0,0){\includegraphics[width=6in]{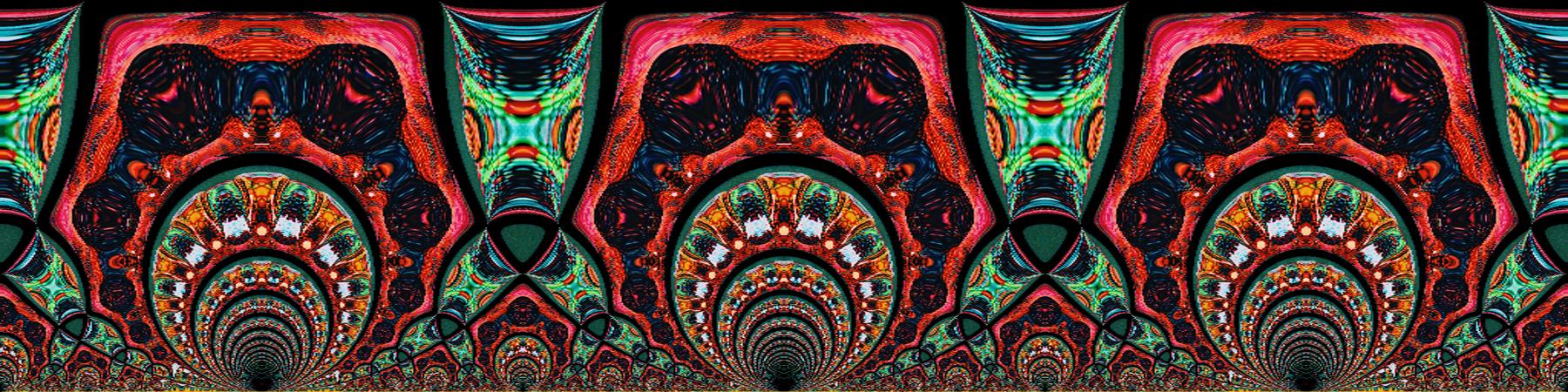}}
\end{picture}}}
\end{picture}
\caption{\small Top: \emph{Blugold Fireworks} design, using \emph{Buoy} as color map. Middle and Bottom: Two designs using \emph{Reptile}. \label{fig:HyperbolicDesigns}}
\end{figure}

We have also created designs using the symmetries in the geometry of the hyperbolic upper half-plane. Some
of the designs we have created are shown in \Fig~\ref{fig:HyperbolicDesigns}. The symmetries in these designs
are much different than those we discussed above for the Euclidean plane. The design at the top of \Fig~\ref{fig:HyperbolicDesigns} is entitled \emph{Blugold Fireworks}. It was exhibited as part
of the 2018~\textsc{Mathematical Art Exhibition} held in San Diego~\cite{MathArtExhibit2018}.
\noindent
The other two designs are more recent creations. One quite interesting feature of these designs, from a mathematical perspective,
is that they display some of the principal geometric objects in the hyperbolic upper half-plane. We will describe what
these principal objects are, and how the designs are constructed. But in order to do so, we first recount
the basic mathematics underlying the geometry of the hyperbolic upper half-plane. References for
additional details are \cite{ref:ClimKotek,Farris2015,ref:Katok,Needham,Terras}.
\subsection{Geometry of the Hyperbolic Upper Half-Plane\label{subsec:GeometryUpperHalfPlane}}
The hyperbolic upper half-plane, $\H$, is the subset of $\C$ defined as follows:
\begin{align}
\H &= \{x + iy\, |\, x\in\R, y > 0\}\label{eq:DefinitionOfHyperbolicHalfPlaneA}\\
\text{with differential metric}\notag\\
ds &= \frac{\sqrt{dx^2 + dy^2}}{y}\label{eq:DefinitionOfHyperbolicHalfPlaneB}
\end{align}
With this metric, the \Emph{length} $\ell(\gamma)$ of a smooth curve $\gamma (t) = x(t) + iy(t)$, $a\leq t \leq b$, is defined as
\[
\ell(\gamma) = \int^b_a \frac{\sqrt{x'(t)^2 + y'(t)^2}}{y(t)}\, dt \; = \int^b_a \frac{|z'(t)|}{\Im{z(t)}}\, dt
\]
This metric is related to the metric $ds_{\Esub}$ for the Euclidean plane, restricted to $\EEL = \{(x,y)\, |\, x\in\R, y > 0\}$, defined by
\begin{equation}\label{eq:EuclideanMetric}
ds_{\Esub} = \sqrt{dx^2 + dy^2}
\end{equation}
In other words, $ds = ds_{\Esub} / y$. This relation is crucial to verifying a number of important facts about the geometry of $\H$.  We begin by discussing the isometries of $\H$.
\subsubsection{Isometries of $\H$\label{subsubsec:IsometriesOfH}}
The isometries of $\H$ are mappings $f\colon \H\to\H$ that preserve the differential metric $ds$. We will show in the next theorem that the set
\begin{equation}\label{eq:IsometriesOfH}
\MG = \left\{f(z) = \frac{az+b}{cz+d}\;\colon\;a, b, c, d\in \R\ \text{with}\ ad - bc = 1\right\}
\end{equation}
contains all the holomorphic isometries of $\H$.

There is also a notion of area in $\H$. The \Emph{area} of a region $\mathcal{U}$ will be denoted by $A(\mathcal{U})$. The area differential $dA$ in $\H$ is given by
\begin{equation}\label{eq:AreaDifferentialForH}
dA = \frac{1}{y^2}\, dx\, dy
\end{equation}
and so we compute $A(\mathcal{U})$ by
\begin{equation}\label{eq:AreaIntegralForH}
A(\mathcal{U}) = \int_{\mathcal{U}} \frac{1}{y^2}\, dx\, dy
\end{equation}
for a suitable region $\mathcal{U} \subset \H$. By a suitable region, we mean any region $\mathcal{U}\subset \H$ for which the integral in \eqref{eq:AreaIntegralForH} is defined, say, as a Riemann integral. The area differential $dA$ makes sense, by a dimensional argument, when viewed
as $dA = dA_{\Esub} / y^2$ and noting that the length differential satisfies $ds = ds_{\Esub} / y$. It also follows from basic facts of Riemannian geometry:
\[
\begin{matrix}ds^2 = \displaystyle\sum_{i,j = 1}^n g_{i,j} dx^i dx^j \\[5pt]\hspace{-0pt}\text{\footnotesize (differential metric squared)}\end{matrix}
\quad
\begin{matrix}\implies \\\phantom{1}\end{matrix}
\quad
\begin{matrix}d\Omega = \displaystyle\sum_{i,j=1}^n \sqrt{\det (g_{i,j})}\; dx^1 \wedge \dots \wedge dx^n \\[5pt]\hspace{-69pt}\text{\footnotesize (differential volume element)}\end{matrix}
\]
as shown in \cite[p.~188, 241]{ref:Boothby}. In this case, $(g_{i,j}) = \begin{pmatrix}1/y^2 & 0 \\ 0 & 1 / y^2\end{pmatrix}$, and we obtain
$ds$ and $dA$ as defined in \eqref{eq:DefinitionOfHyperbolicHalfPlaneB} and \eqref{eq:AreaDifferentialForH}. We will show that isometries in $\MG$
also preserve area in $\H$.
\begin{theorem}\label{thm:IsometriesInH}\textsl{The set $\MG$ defined in~\eqref{eq:IsometriesOfH} contains all the holomorphic isometries of $\H$. In addition to preserving the metric differential $ds$, these isometries also preserve the area differential $dA$. The isometries in $\S_\H$ are also described by
\begin{equation}\label{eq:IsometriesOfHAlternative}
\MG = \left\{f(z) = \frac{az+b}{cz+d}\;\colon\;a, b, c, d\in \R\ \text{with}\ ad - bc > 0\right\}
\end{equation}}
\end{theorem}
\begin{small}\begin{proof}
For $f(z) = \dfrac{az+b}{cz+d}$ in $\MG$, we have
\begin{align*}
2i\,\Im{f(z)} &= \frac{(ad-bc)(z-\overline{z})}{|cz+d|^2}\\
&= \frac{1}{|cz+d|^2}\, 2i\,\Im{z}
\end{align*}
and therefore $\Im{f(z)} = \Im{z}\, /\, |cz+d|^2$. Consequently, $f(z)\in\H$ if and only if $z\in\H$. (Note: $cz + d = 0$ is
only possible when $z = -d/c \in\R$ and such $z$ are not in $\H$.)

Now, for $f(z) = u + iv$ with $u\in\R$ and $v > 0$, we have $v = y / |cz+d|^2$. We also have
\begin{align*}
du^2 + dv^2 &= |J|(dx^2 + dy^2)\\
&= \left[\left(\frac{\partial u}{\partial x}\right)^2 + \left(\frac{\partial u}{\partial x}\right)^2\right]
(dx^2 + dy^2) \\
&= |f'(z)|^2 \left(dx^2 + dy^2\right)\\
&= \frac{1}{|cz+d|^4}\, \left(dx^2 + dy^2\right)
\end{align*}
where we made use of the Cauchy-Riemann equations to simplify the Jacobian $|J|$ for the change of variables in the
first line. Consequently, we obtain
\begin{align*}
\frac{du^2 + dv^2}{v^2} &= \frac{|cz+d|^4}{y^2}\,\cdot\,\frac{1}{|cz+d|^4}\, \left(dx^2 + dy^2\right)\\
&= \frac{dx^2+dy^2}{y^2}
\end{align*}
and that shows that $f(z)$ is an isometry of $\H$. Moreover, in our calculations we computed $f'(z)$ for $z\in\H$,
so $f$ is holomorphic on $\H$.

To prove preservation of $dA$, we calculate as above:
\begin{align*}
du\, dv = &= |J|\, dx\, dy\\
 &= \frac{1}{|cz+d|^4}\, dx\, dy
\end{align*}
and $du\, dv /v^2 = dx\, dy / y^2$ follows just as above. Thus, the isometry $f$ also preserves the area differential $dA$.

Now, suppose that $f(z) = \dfrac{az + b}{cz+d}$ with $ad-bc > 0$. Let $r^2 = ad-bc$. Then,
\begin{align*}
f(z) &= \frac{r^2}{r^2}\,\frac{a' z + b'}{c' z + d'}\\
&= \frac{a' z + b'}{c' z + d'}
\end{align*}
with $a'd' - b'c' = 1$. Therefore, $f(z)$ is in $\MG$, as defined in \eqref{eq:IsometriesOfH}. Since the
reverse inclusion obviously holds, it follows that $\MG$ is described by both~\eqref{eq:IsometriesOfH}
and \eqref{eq:IsometriesOfHAlternative}.

Finally, by \cite[Theorem~5]{ref:Cima}, all holomorphic mappings
from the disc $\D = \{z\;\colon\, |z| < 1\}$ to itself have the form $F(z) = t(z - c)/(1-\overline{c}z)$, for some
$t,c\in \C$ with $|t| = 1$ and $|c|<1$. Conjugating each $F$ with the conformal map $g\colon \D \to \H$ given
by $g(z) = \hbox{\Large $\frac{z+i}{iz + 1}$}$, we obtain all the functions $f = g\circ F\circ g^{-1}$ that belong to $S_\H$.
We omit the details for verifying this last statement, because we will not be using the fact that
$\MG$ consists of \Emph{all} the holomorphic isometries of $\H$. Complete details are in
\cite[Theorem~2.4, p.~222]{ref:Stein}.
\end{proof}\end{small}

\begin{small}\begin{remark}\label{rem:OtherIsometries}
Theorem~\ref{thm:IsometriesInH} deals with the holomorphic isometries of $\H$. There are other isometries.
For example, the function $f(x+iy) = -x + iy$ is an isometry, since it clearly preserves $ds$. However,
it is not holomorphic on $\H$ due to its failure to satisfy the Cauchy-Riemann equations.
\end{remark}\end{small}
\par\noindent
Theorem~\ref{thm:IsometriesInH} tells us that these sets of mappings are all isometries:
\begin{enumerate}
\item $\M_\rho\colon z \to \rho z$, for $\rho > 0$. In $\EEL$\!, this mapping would be a similarity transformation when
$\rho \neq 1$, not an isometry. In $\H$, however, this mapping is an isometry for all $\rho > 0$.

\item $\tau_u\colon z \to z + u$, for $u\in\R$. Thus, all horizontally oriented translations are isometries of $\H$.

\item $\I_r\colon z \to \dfrac{-r}{z/r}$ for $r > 0$. This operation is called \emph{inversion} through the
circle of radius $r$, center $0$ in $\C$. However, it is also inversion through the upper semicircle
in $\H$ defined by $\Sur{r}{0} = \{z\in\H\;\colon\, |z| = r\}$. In other words,
$\Sur{r}{0} = \{x + iy \;\colon\, x^2 + y^2 = r, y > 0\}$. In subsequent work, we shall also deal with upper semicircles of radius $r$ and center $u\in\R$, which we denote by $\Sur{r}{u}$. Note that $\Sur{r}{u} = \tau_u\left(\Sur{r}{0}\right)$ and $\Sur{r}{0} = \M_r\left(\Sur{1}{0}\right)$.
\end{enumerate}
These special isometries generate all the isometries in $\MG$ through composition. In fact, if $c \neq 0$, then
by long division we obtain
\begin{align*}
\frac{az+b}{cz+d} &= \frac{a}{c} + \frac{b - da/c}{cz + d}\\
&= \frac{a}{c} + \frac{1}{c^2}\,\frac{-1}{z + d/c}
\end{align*}
hence $f = \tau_{a/c}\circ \M_{1/c^2}\circ \I_1 \circ \tau_{d/c}$. While if $c = 0$, then
\[
\frac{az+b}{d} = \frac{a}{d}\, z + \frac{b}{d}
\]
Since $c = 0$, $ad - bc = 1$ reduces to $ad = 1$, and we have
\[
f(z) = a^2 (z + b/a)
\]
Thus, $f = M_{a^2}\circ \tau_{b/a}$.
\subsubsection{Geodesics\label{subsubsec:Geodesics}}
A \Emph{geodesic} in $\H$, connecting two points $z_1$ and $z_2$, is a (piecewise) smooth curve $\gamma (t) = x(t) + iy(t)$
for $a \leq t \leq b$ that satisfies $\gamma (a) = z_1$, $\gamma(b) = z_2$, and which has minimum length $\ell(\gamma)$.
We shall now prove that, in $\EEL$, these geodesics lie along vertical rays or semicircles.

\begin{theorem}\textsl{%
The geodesics in $\H$ lie along the following two types of curves in $\EEL$:
\textbf{(1)} vertical rays $\Ru{u}$ emanating from $\R$:
$\Ru{u} = \{u + iy\;\colon\; y > 0,\ \text{fixed}\ u\in\R\}$, or \textbf{(2)} open semicircles $\Sur{r}{u}$ centered
on $\R$: $\Sur{r}{u} = \{x + iy\;\colon\; (x-u)^2 + y^2 = r^2, y > 0,\ \text{fixed}\ u\in\R, r > 0\}$.}
\end{theorem}
\begin{small}\begin{proof}
First, we consider $z_1 = iy_1$ and $z_2 = iy_2$, choosing subscripts so that $y_2 > y_1 > 0$. For a smooth curve
$\gamma(t) = x(t) + iy(t)$ satisfying $\gamma(a) = y_1$ and $\gamma(b) = y_2$, we have
\begin{align*}
\ell(\gamma) &= \int^b_a \frac{\sqrt{x'(t)^2 + y'(t)^2}}{y(t)}\, dt\\
&\geq \int^b_a \frac{y'(t)}{y(t)}\, dt\\
&= \ln(y_2/y_1)
\end{align*}
Moreover, this lower bound of $\ln(y_2/y_1)$ is realized for $\gamma(t) = [y_1 + (y_2-y_1)t]i$ for $0\leq t \leq 1$.
Therefore, this function $\gamma$ is a geodesic in $\H$, and clearly it lies on the ray $\Ru{0}$. The minimum
property also extends to the class of all continuous, piecewise smooth curves, by splitting integrals over $[a,b]$
into finite sums of integrals.  
Since the horizontal translation $\tau_u$ is an isometry, conjugation with $\tau_u$ implies that
geodesics also lie along each vertical ray $\Ru{u}$.

Second, we consider two points $z_1 \neq z_2$ on the open semicircle $\Sur{1}{0}$, having
$\arg(z_2) = \theta_2 > \theta_1 = \arg(z_1)$. The isometry $f(z) = (z + 1)/(-z + 1)$ maps
$\Sur{1}{0}$ to $\Ru{0}$, with $iy_2 = f(z_2)$, $iy_1 = f(z_1)$, and $y_2 > y_1$.
Given a geodesic $\gamma(t) = i[(y_1 + (y_2 - y_1)t]$ for $0\leq t \leq 1$, connecting $iy_1$ and $iy_2$
on $\Ru{0}$, we apply the isometry $f^{-1}$ to obtain $f^{-1}\circ \gamma$ as a geodesic
on $\Sur{1}{0}$ connecting $z_1$ and $z_2$. Thus, $\Sur{1}{0}$ contains geodesics in $\H$. Since
the isometry $\M_r$ maps $\Sur{1}{0}$ to $\Sur{r}{0}$ it follows that $\Sur{r}{0}$ contains geodesics in $\H$.
Finally, since the isometry $\tau_u$ maps $\Sur{r}{0}$ to $\Sur{r}{u}$, it follows that $\Sur{r}{u}$ contains geodesics in $\H$.
\end{proof}\end{small}

\begin{small}\begin{remark}\label{rem:DefnOfHDistance}
The \Emph{distance} $d(z_1, z_1)$ between two points $z_1, z_2 \in\H$ is defined to be the length of a
geodesic that connects $z_1$ and $z_2$. For example, we found above that $d(x+iy_1, x+iy_2) = \ln (y_2/y_1)$
for $y_2 > y_1 > 0$. In general, for $x + iy_1, x + iy_2 \in \H$, we have $d(x+iy_1, x+iy_2) = |\ln (y_2/y_1)|$.
It is important to note that $|\ln (y_2/y_1)| \to \infty$ if either $y_1 \to 0$ or $y_2\to 0$. Consequently, the real line $\R$ is a \emph{line at infinity} for all points in $\H$.

There is a distance formula for all $z, w \in\H$, given by
\begin{equation}\label{eq:GeneralDistanceFormulaInH}
d(z, w) = \ln \frac{|z - \overline{w}| + |z - w|}{|z - \overline{w}| - |z - w|}
\end{equation}
but we will not need it. Interested readers will find a proof of \eqref{eq:GeneralDistanceFormulaInH}
in Katok~\cite[Theorem~1.2.6, p.~6]{ref:Katok}.
\end{remark}\end{small}

\begin{figure}[!htb]\centering
\setlength{\unitlength}{1in}
\begin{picture}(5.63,3.6)
\put(0,1.9){\includegraphics[width=5.63in]{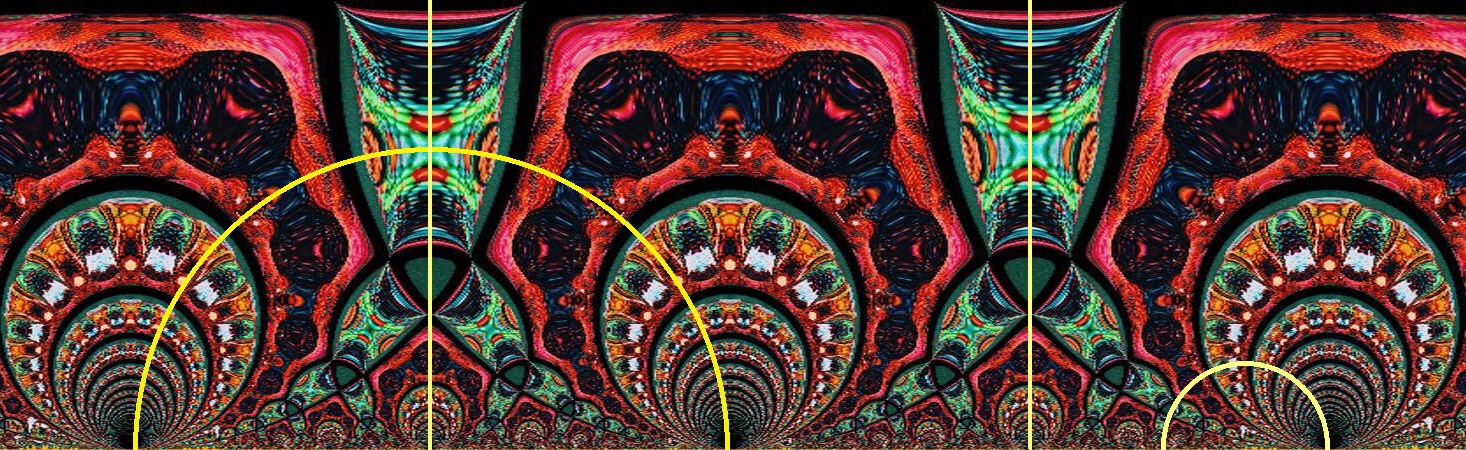}}
\put(0,0){\includegraphics[width=1.8in]{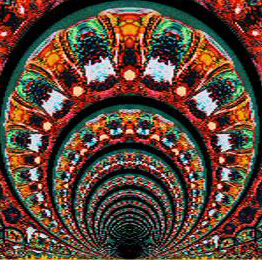}}
\put(1.9,0){\includegraphics[width=1.8in]{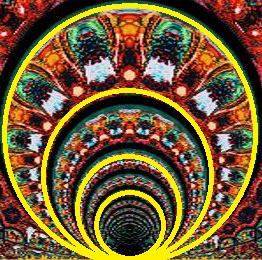}}
\put(3.8,0){\includegraphics[width=1.83in]{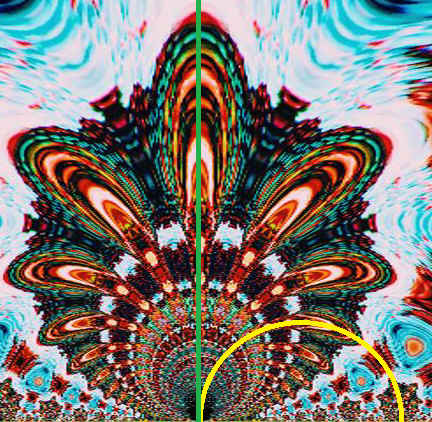}}
\end{picture}
\caption{\small Top: Four geodesics, drawn in yellow over one of our hyperbolic designs.
The two intersecting geodesics on the left are of types $\Ru{u}$ and $\Sur{r}{u}$. The two
disjoint geodesics on the right are also of types $\Ru{u}$ and $\Sur{r}{u}$. Bottom Left: Portion of the
same hyperbolic design that contains horocycles. Bottom Middle: Yellow circles indicating some of these
horocycles. These horocycles are orthogonal to geodesics. Bottom Right: Yellow semicircle
and green vertical line indicating geodesics on the middle design in \Fig~\ref{fig:HyperbolicDesigns}. \label{fig:HyperbolicDesignWithHorocycles}}
\end{figure}

For simplicity, in the rest of the paper, we shall refer to rays of type $\Ru{u}$ and
open semicircles of type $\Sur{r}{u}$ as \emph{geodesics}. Strictly speaking, they contain geodesics, but
there is little chance of confusion and our language is more straightforward if we simply call them
geodesics as well. These geodesics in $\H$ can be interpreted as a model for the undefined term
\emph{lines} referred to in postulates of geometry. In fact, these geodesics do satisfy the first four
of Euclid's postulates. However, they violate the notion of Euclidean parallelism. For example, in
the image at the top of \Fig~\ref{fig:HyperbolicDesignWithHorocycles}, the two geodesics on the left intersect
at a point. Yet, they fail to intersect the vertical geodesic on the right of the image. This situation violates
the uniqueness of a parallel line, through a point not on a line, required in Euclidean geometry.

Returning to the artworks in \Fig~\ref{fig:HyperbolicDesigns}, it is interesting
that parts of these designs correspond to geodesics. On the top of \Fig~\ref{fig:HyperbolicDesignWithHorocycles} we have shown
that geodesics of both types, $\Sur{r}{u}$ and $\Ru{u}$, are evident within the design shown at the bottom of \Fig~\ref{fig:HyperbolicDesigns}. Parts of the other two designs also correspond to these types of geodesics. For instance, on the
bottom right image in \Fig~\ref{fig:HyperbolicDesignWithHorocycles} we have shown how a part of the
middle design in \Fig~\ref{fig:HyperbolicDesigns} corresponds to both types of geodesics in $\H$.
\par
The images at the bottom of \Fig~\ref{fig:HyperbolicDesignWithHorocycles} also contain geometric objects
related to geodesics in $\H$. These objects are circles that are tangent to $\R$ at one
point and have all other points lying in $\H$. To be specific, for $u\in\R$ and $r > 0$, a \Emph{horocycle} $H_r(u)$
is defined by
\[
H_r(u) = \{x + iy\;\colon\; (x-u)^2 + (y-r)^2 = r^2, y > 0\}
\]
so in $\EEL$ it is a circle with center $(u,r)\in\R^2$, and radius $r$, but omitting the point $(u,0)$ on the $x$-axis.
These horocycles are \Emph{not} geodesics. However, we will now discuss how the family of all horocycles are
orthogonal curves in $\H$ for the family of all geodesics.
\subsubsection{Angles and Conformality in $\H$, Horocycles and Geodesics}\label{subsubsec:AnglesAndConformalityInH}
An \Emph{angle} in $\H$ is defined to be an angle between tangent vectors of two
curves meeting at a point. The following theorem shows that these angles are the same in both $\H$ and $\EEL$.
\begin{theorem}\label{thm:AnglesSameinHandE}\textsl{%
Let $\theta_\H$ and $\theta_{\Esub}$ be the angles between two curves at some intersection point in
$\H$ and $\EEL$, respectively. Then, $\theta_\H = \theta_{\Esub}$.}
\end{theorem}
\begin{small}\begin{proof}
We can write the infinitesimal quadratic form $ds^2 = (dx^2 + dy^2)/y^2$ as
\begin{align*}
ds^2 &= \frac{\bigl\langle [dx, dy], [dx, dy]\bigr\rangle_{\Esub}}{y^2}\\
&=\|[dx, dy]\|^2_\H
\end{align*}
where $\bigl\langle [dx, dy], [dx, dy]\bigr\rangle_{\Esub}$ stands for the standard inner product of the vector of
differentials $[dx, dy]$ with itself, and $\|[dx,dy]\|^2_\H$ is our notation for $ds^2$ thought of as a
quadratic form of the vector $[dx, dy]$. The inner product $\bigl\langle [dx, dy], [d\tilde{x}, d\tilde{y}]\bigr\rangle_\H$ corresponding to this quadratic form is then
\[
\bigl\langle [dx, dy], [d\tilde{x}, d\tilde{y}]\bigr\rangle_\H =
\frac{dx\, d\tilde{x} + dy\, d\tilde{y}}{y^2}
\]
Consequently, $\cos \theta_\H$ satisfies
\begin{align*}
\cos \theta_\H &= \frac{\bigl\langle [dx, dy], [d\tilde{x}, d\tilde{y}]\bigr\rangle_\H}
{\|[dx,dy]\|_\H\, \|[d\tilde{x},d\tilde{y}]\|_\H}\\
&= \frac{dx\, d\tilde{x} + dy\, d\tilde{y}}{\sqrt{dx^2 + dy^2}\,\sqrt{d\tilde{x}^2 + d\tilde{y}^2}}\\
&= \cos \theta_{\Esub}
\end{align*}
Thus, we must have $\theta_\H = \theta_\EEL$.
\end{proof}\end{small}
Since angles in $\H$ and $\EEL$ always correspond to angles between
tangent vectors, we have proved that angles in the two geometries are
always the same. The two geometries are said to be \Emph{conformal}.
\begin{corollary}\label{cor:IsometriesForHPreserveAngles}
The isometries in $\MG$ preserve angles in $\H$.
\end{corollary}
\begin{small}\begin{proof}
If $f(z) = (az+b)\,/\,(cz+d) \in\MG$, then $f'(z) = 1/|cz+d|^2 \neq 0$. Therefore,
$f$ is a conformal mapping on $\EEL$ by \cite[Theorem~3]{ref:Cima}. Hence Theorem~\ref{thm:AnglesSameinHandE}
above implies that $f$ is a conformal mapping on $\H$.
\end{proof}\end{small}
\par\noindent
This corollary can also be proved using the identity
\[
2\bigl\langle \mathbf{v}, \mathbf{w}\bigr\rangle = \|\mathbf{v} + \mathbf{v}\|^2
- \|\mathbf{v}\|^2 - \|\mathbf{w}\|^2
\]
relating inner products and quadratic forms, and the fact that an isometry preserves the
quadratic form $ds^2$. However, our proof highlights the relation between the geometries
of $\EEL$ and $\H$.

We now return to the concept of horocycles in $\H$, and how they are illustrated in the designs
shown in \Fig~\ref{fig:HyperbolicDesigns}. The simplest type of horocycles are the sets of
form, $\{t + iv\;\colon\; t\in\R\}$, parameterized by varying $iv$ with $v>0$. These sets
are horizontal lines in $\EEL$, but in $\H$ they are \Emph{not} geodesics.
Each geodesic ray $\Ru{u}$, for $u\in\R$, lies orthogonal in $\EEL$ at each of its points to one of
these horizontal lines, and therefore each geodesic ray $\Ru{u}$ also lies orthogonal in $\H$ at each
of its points to one of these horocycles $\{t + iv, t\in\R\}$.

The second type of horocycles are those that lie orthogonal to points of open semicircle
geodesics. The isometry $\I_r (z) = -r^2 / z$ maps the geodesic ray $\Ru{0}$ to itself (with $ri$ held fixed),
and maps the geodesic ray $\Ru{-r}$ to the open semicircular geodesic $\Sur{r/2}{r/2}$. The horocycles
for $\Ru{-r}$, expressed as $\{t + ryi\;\colon\; t\in\R\}$ for each $y > 0$, are mapped by
$\I_r$ to sets of the form $\{w\in\H\;\colon\; \left|w - \frac{r}{2y}\, i\right|^2 = \left(\frac{r}{2y}\right)^2\}$, which
are circles in $\EEL$ except for the one point $0 + 0i \notin\H$. By Corollary~\ref{cor:IsometriesForHPreserveAngles},
these horocycles are orthogonal to the open semicircular geodesic $\Sur{r/2}{r/2}$ at all of its points.
They are circles in $\EEL$ that are tangent to the point $(0,0)$, and all their points excepting $(0,0)$ lie in $\H$. Conjugating with horizontal translation $\tau_u$ for any fixed $u\in\R$, we find that the horocycles for all open semicircular geodesics in $\H$ are circles in $\EEL$ except for one point that is tangent to $\R$. Since these horocycles are all tangent
to $\R$ at the same point $u$, with all radii $r > 0$, it follows that each family of horocycles is also orthogonal to
the geodesic ray $\Ru{u}$ at each of its points. On the bottom left and bottom middle of \Fig~\ref{fig:HyperbolicDesignWithHorocycles}, we illustrate a collection of such horocycles in one of our designs. The geodesics  drawn on the design at the top of this figure are orthogonal at each of their points to such horocycles. On the bottom right of \Fig~\ref{fig:HyperbolicDesignWithHorocycles}, we show a part
of the middle design in \Fig~\ref{fig:HyperbolicDesigns} that exhibits both horocycles and geodesics.
\subsection*{Creating designs with hyperbolic symmetry\label{sec:CreatingDesignsHyperbolicSymmetry}}
Our method for creating designs with hyperbolic symmetries is similar to our method for Euclidean symmetries. We symmetrize
a given function $f$ with domain $\H$. The symmetries will be a subgroup of $\MG$. We cannot use $\MG$ itself because
the only functions on $\H$, symmetric with respect to all the transformations in $\MG$, are constant functions.
Following Farris~\cite{Farris2015}, we will use the subgroup $\MGZ$ known as the \emph{modular group}. The
modular group $\MGZ$ is defined as
\begin{equation}\label{eq:DefnOfPSL2Z}
\MGZ = \left\{f(z) = \frac{jz + k}{mz+n}\;\colon\;j, k, m, n \in\Z, jn - mk = 1\right\}
\end{equation}
Note that $jn - mk$ is the determinant of the matrix $\begin{pmatrix}j & k\\ m & n\end{pmatrix}$ of coefficients of
$f(z) = \dfrac{jz + k}{mz+n}$. The set $\MGZ$ is a group because composition of two members $f(z) = \dfrac{jz+k}{mz+n}$ and
$g(z) = \dfrac{j' z + k'}{m' z + n'}$ satisfies
\[
(f\circ g) (z) = \frac{(j j' + k m')z + (j k' + k n')}{(m j'+ n m')z + (m k' + n n')}
\]
which corresponds to multiplication of the matrices of coefficients of $f$ and $g$:
\[
\begin{pmatrix}
j & k\\ m & n
\end{pmatrix}
\begin{pmatrix}
j' & k'\\ m' & n'
\end{pmatrix}
=
\begin{pmatrix}
jj' + k m' & j k' + k n'\\ m j' + n m' & m k' + n n'
\end{pmatrix}
\]
and we know that determinants of matrices respect multiplication. So the determinants of each of the matrices
in the equation above satisfy $1 \cdot 1 = 1$, hence $f\circ g$ is a member of $\MGZ$. Furthermore,
$f^{-1}(z) = \dfrac{nz - k}{-mz + j}$, and therefore $f^{-1}\in\MGZ$.
\begin{small}\begin{remark}\label{rem:GammaAsMatrixGroup}
The group $\MG$ is isomorphic to a subgroup of the matrix factor group:
\[
\mathrm{PSL}(2,\R) = \mathrm{SL}(2,\R)\, /\, \{\mathrm{Id}, -\mathrm{Id}\}
\]
where $\mathrm{SL}(2,\R)$ is the \Emph{special linear group} consisting of all $2$~by~$2$ matrices over $\R$ having determinant $1$, and $\mathrm{Id}$ is the $2$~by~$2$ identity matrix. The group $\mathrm{PSL}(2,\R)$ is related to the projective geometry
of all lines through the origin~\cite[p.~179]{ref:ClimKotek}. It is called the
\Emph{projective special linear group} over $\R$. The group $\MGZ$ is isomorphic to
$\mathrm{PSL}(2,\Z) = \mathrm{SL}(2,\Z) \,/\, \{\mathrm{Id}, -\mathrm{Id}\}$,
where $\mathrm{SL}(2,\Z)$ consists of matrices in $\mathrm{SL}(2,\R)$ with integer coefficients.
\end{remark}\end{small}
\par\noindent
The special isometries, $\M_r$, $\tau_u$, $\I_r$, mostly are not members of $\MGZ$ due to the requirement that
their coefficients belong to $\Z$. In fact, the special isometries that belong to $\MGZ$ are
\begin{enumerate}
\item Translations: $T^n\colon z \to z + n$, for $n\in\Z$. The unit-translation $T^1$ will be written as just $T$.
These translations obey the group operation in $\MGZ$: $T^m \circ T^n = T^{m+n}$ for
all $m, n \in \Z$. 
\item Inversion: $\I\colon z \to -1/z$ 
\end{enumerate}
Compositions with these isometries are sufficient to generate all the isometries in $\MGZ$.
\begin{theorem}\label{thm:IandTgenerateGamma}\textsl{%
The unit translation $T$ and inversion $\I$ generate $\MGZ$. More precisely, if $f\in\MGZ$, then
$f$ can be written as some interlaced composition of $\I$ with various translations $T^p$:
\begin{equation}\label{eq:DecompositionWithTandI}
f =  \left(\prod_{\ell=1}^L T^{p_\ell} \I\right) T^{p_0}
\end{equation}}
\end{theorem}
\begin{small}\begin{proof}
Let $f(z) = \hbox{\normalsize $\frac{jz+k}{mz+n}$}$ be an arbitrary function in $\MGZ$. First, suppose $m = 0$. Then $f(z) = (j/n)z + (k/n)$ and we have $j n = 1$. Consequently,
$(j, n) = (1,1)$ or $(-1,-1)$. Hence, either $f(z) = z + k$ or $f(z) = z - k$, and so $f = T^k$ or $f = T^{-k}$. Now, if $m \neq 0$, we reduce
to the first case as follows. We have, where $\ell \in\Z$,
\[
\I f (z) = \frac{-mz - n}{jz + k}\, , \qquad T^\ell f (z) = \frac{(j+\ell m)z + (k + \ell n)}{mz + n}
\]
If $|j| > |m|$, the division algorithm gives $j = mq + r$ with $0 \leq |r| < |m|$. Hence $T^{-q} f(z) = \hbox{\normalsize $\frac{rz + (k + \ell n)}{mz + n}$}$. Then apply $\I$
to obtain $\I T^{-q} f (z) = \hbox{\normalsize $\frac{-mz - n}{rz + (k + \ell n)}$}$ with $|r| < |m|$. Applying powers of $T$, followed by $\I$, eventually results in a remainder $r = 0$ as coefficient of $z$ in the denominator. That is, we arrive at $T^{p_0}$ for some $p_0\in\Z$. Thus, we obtain
\[
\left(\prod^L_{s = 1} \I T^{-q_s}\right) f = T^{p_0}
\]
Solving for $f$, we obtain the result in Equation~\eqref{eq:DecompositionWithTandI}.
\end{proof}\end{small}
The symmetrization $f_\S$ of a function $f$, defined on $\H$, can be done as follows
\begin{equation}\label{eq:HyperbolicSymmetrizationFunction}
f_\S (z) = \sum_{g\in \MGZ} f\bigl(g(z)\bigr)
\end{equation}
Since $\MGZ$ is a group, the function $f_\S$ is guaranteed to satisfy the symmetry condition $f(g(z)) = f(z)$ for all $g\in \MGZ$. In practice, of course, we can typically only create partial sums of the infinite series for $f_\S$. Nevertheless, as shown above, our designs using such partial sums display many important features of the geometry of $\H$.

To create designs by domain coloring of $f_\S(z)$, we need to express the series defining $f_S$ in a more convenient form. To do that, we observe that the condition $jn - mk = 1$ can be rewritten as
\begin{equation}\label{eq:RelativelyPrimeCondition}
nj + (-m)k = 1
\end{equation}
Equation~\eqref{eq:RelativelyPrimeCondition} is a famous one from Number Theory. It is equivalent to the integers $j$ and $k$ being \emph{relatively prime}, i.e., their greatest common divisor is $1$, which we write as $\mathrm{gcd}(j, k) = 1$. The numbers $n$ and $-m$ are called \emph{B\'ezout coefficients} for $j$ and $k$,
and they ensure that Equation~\eqref{eq:RelativelyPrimeCondition} holds. Because
\[
\mathrm{det}\begin{pmatrix}j & k\\ m & n\end{pmatrix} = 1
\]
will hold if $j$ and $k$ are perturbed by $j = \ell m$ and $k = \ell n$ for $\ell\in\Z$, these B\'ezout coefficients
determine all transformations $g$ having the form
\[
g(z) = \ell + \frac{jz + k}{mz+n}
\]
where $(n, -m)$ is a single pair of B\'ezout coefficients for $(j,k)$. Therefore, we will assume that the initial
function $f$ has period $1$ in the $x$-variable, and write $f_\S$ as
\begin{equation}\label{eq:HyperbolicSymmetrizationFunctionFormula2}
f_\S (z) = \sum_{\mathrm{gcd}(j,k) = 1} f\left(\frac{jz + k}{mz+n}\right)
\end{equation}

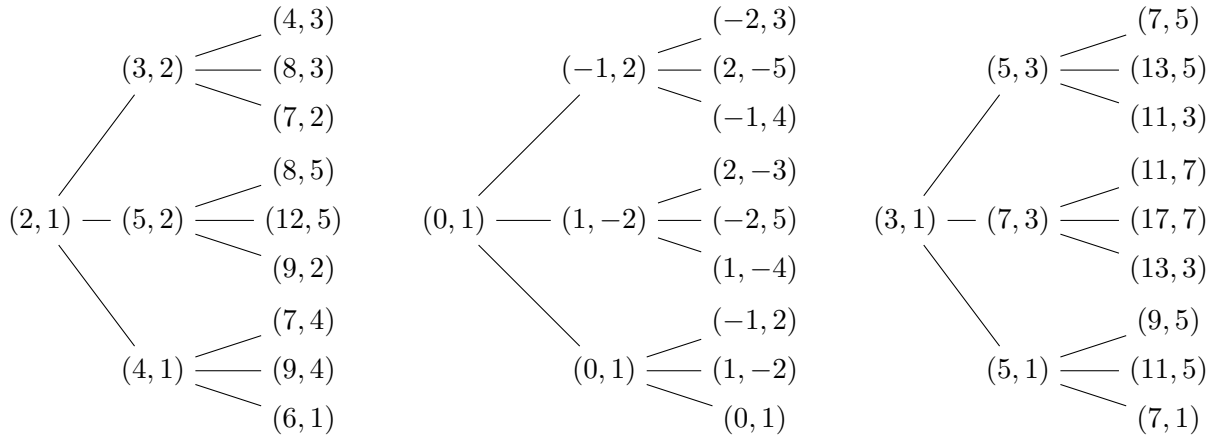
\begin{figure}[!htb]\centering
\begin{tikzpicture}
\node at (-7.5,0) (A){$(2,1)$};
\node at (-6,2) (B_1){$(3,2)$};
\node at (-6,0) (B_2){$(5,2)$};
\node at (-6,-2) (B_3){$(4,1)$};
\node at (-4,2.65) (C_1){$(4,3)$};
\node at (-4,2) (C_2){$(8,3)$};
\node at (-4,1.35) (C_3){$(7,2)$};
\node at (-4,.65) (C_4){$(8,5)$};
\node at (-4,0) (C_5){$(12,5)$};
\node at (-4,-.65) (C_6){$(9,2)$};
\node at (-4,-1.35) (C_7){$(7,4)$};
\node at (-4,-2) (C_8){$(9,4)$};
\node at (-4,-2.65) (C_9){$(6,1)$};
\draw (A) -- (B_1);
\draw (A) -- (B_2);
\draw (A) -- (B_3);
\draw (B_1) -- (C_1);
\draw (B_1) -- (C_2);
\draw (B_1) -- (C_3);
\draw (B_2) -- (C_4);
\draw (B_2) -- (C_5);
\draw (B_2) -- (C_6);
\draw (B_3) -- (C_7);
\draw (B_3) -- (C_8);
\draw (B_3) -- (C_9);
\node at (-2,0) (X){$(0,1)$};
\node at (0,2) (Y_1){$(-1,2)$};
\node at (0,0) (Y_2){$(1,-2)$};
\node at (0,-2) (Y_3){$(0,1)$};
\node at (2,2.65) (Z_1){$(-2,3)$};
\node at (2,2) (Z_2){$(2,-5)$};
\node at (2,1.35) (Z_3){$(-1,4)$};
\node at (2,.65) (Z_4){$(2,-3)$};
\node at (2,0) (Z_5){$(-2,5)$};
\node at (2,-.65) (Z_6){$(1,-4)$};
\node at (2,-1.35) (Z_7){$(-1,2)$};
\node at (2,-2) (Z_8){$(1,-2)$};
\node at (2,-2.65) (Z_9){$(0,1)$};
\draw (X) -- (Y_1);
\draw (X) -- (Y_2);
\draw (X) -- (Y_3);
\draw (Y_1) -- (Z_1);
\draw (Y_1) -- (Z_2);
\draw (Y_1) -- (Z_3);
\draw (Y_2) -- (Z_4);
\draw (Y_2) -- (Z_5);
\draw (Y_2) -- (Z_6);
\draw (Y_3) -- (Z_7);
\draw (Y_3) -- (Z_8);
\draw (Y_3) -- (Z_9);
\node at (4,0) (A){$(3,1)$};
\node at (5.5,2) (B_1){$(5,3)$};
\node at (5.5,0) (B_2){$(7,3)$};
\node at (5.5,-2) (B_3){$(5,1)$};
\node at (7.5,2.65) (C_1){$(7,5)$};
\node at (7.5,2) (C_2){$(13,5)$};
\node at (7.5,1.35) (C_3){$(11,3)$};
\node at (7.5,.65) (C_4){$(11,7)$};
\node at (7.5,0) (C_5){$(17,7)$};
\node at (7.5,-.65) (C_6){$(13,3)$};
\node at (7.5,-1.35) (C_7){$(9,5)$};
\node at (7.5,-2) (C_8){$(11,5)$};
\node at (7.5,-2.65) (C_9){$(7,1)$};
\draw (A) -- (B_1);
\draw (A) -- (B_2);
\draw (A) -- (B_3);
\draw (B_1) -- (C_1);
\draw (B_1) -- (C_2);
\draw (B_1) -- (C_3);
\draw (B_2) -- (C_4);
\draw (B_2) -- (C_5);
\draw (B_2) -- (C_6);
\draw (B_3) -- (C_7);
\draw (B_3) -- (C_8);
\draw (B_3) -- (C_9);
\end{tikzpicture}
\caption{\small Left: Trinary tree of relatively prime positive integers generated by $(2,1)$. Right: Trinary tree
of relatively prime positive integers generated by $(3,1)$. Middle: Trinary tree of B\'ezout coefficients
generated by $(0,1)$, corresponding to both of the other trees. All trees are shown to a depth of $2$.\label{fig:TreesForIsometries}}
\end{figure}

To efficiently calculate the series in Equation~\eqref{eq:HyperbolicSymmetrizationFunctionFormula2}, we use a recursive, tree-based method for computing relatively prime pairs of positive integers and associated pairs of B\'ezout coefficients. Randall~\cite{Randall1994} proves that pairs of relatively prime positive integers can be computed using $F(j,k) = (2j+k, j)$ in the following recursive formulas:
\begin{align*}
F(j,k) &= (2j+k,j)\\
F(k,j) &= (2k+j,k)\\
F(j,-k) &= (2j-k,j)
\end{align*}
starting from either $(2,1)$ or $(3,1)$ as initial pair. This recursive calculation generates two distinct trinary trees with roots $(2,1)$ and $(3,1)$, as illustrated in \Fig~\ref{fig:TreesForIsometries}.
We have found that there is a similar recursive computation for finding associated B\'ezout coefficients that works for \Emph{both} of the nodes $(j,k)$ in these trinary trees. In \cite{ref:BezoutTreePaper}, we show that using $G(u,v) = (v, u - 2v)$ and the recursive equations
\begin{align*}
G(u,v) &= (v, u-2v)\\
G(v,u) &= (u, v - 2u)\\
G(u,-v) &= (-v, u+2v),
\end{align*}
starting from $(u,v)= (0,1)$, generates B\'ezout coefficients $(u,v)$ satisfying $uj + vk = 1$ for each node $(j,k)$ in both trees of relatively prime positive integers. See the tree in the middle of \Fig.~\ref{fig:TreesForIsometries}.
From these B\'ezout coefficients, and their associated pairs $(j,k)$ of relatively prime positive integers, we obtain transformations $g(z) = \frac{jz + k}{mz+n} \in \MGZ$.

Based on this recursive, tree-based organization of elements of $\MGZ$ we compute symmetrized designs
by the following procedure:
\begin{center}
\textbf{Hyperbolic Symmetrized Design Procedure}
\end{center}
\begin{itemize}
\item[1.] Start with a function $f(z)$ having period $1$ in $x$. It corresponds
to $(j,k) = (1,0)$ and $(m, n) = (0, 1)$. Then add $f(-1/z)$, corresponding to $(j,k) = (0,-1)$, and
$(m,n) = (1,0)$. For this step, we have $f(z) + f(-1/z)$.
\item[2.] Run through the trinary trees of relatively prime integers $(j,k)$, starting at the roots
$(2,1)$ and $(3,1)$. For each pair $(j,k)$, and associated B\'ezout coefficients $(m,-n)$ at the corresponding
node in the B\'ezout tree, add the terms $f\left(\frac{jz+k}{mz+n}\right)$ and
$f\left(\frac{-jz+k}{mz-n}\right)$ to the terms already summed.
\item[3.] After adding a large number of terms---we typically used about $400$ terms---create a symmetrized
design using the \Emph{Domain Coloring Procedure} in \Fig~\ref{fig:ColorMapProcedure}.
\end{itemize}
The designs in \Fig~\ref{fig:HyperbolicDesigns} were all created using this method, starting with various
functions $f$. For example, the design in the middle of \Fig~\ref{fig:HyperbolicDesigns} was created using
\[
f(x + iy) = 2i\, y\cos(2\pi x) + 2y\sin(2\pi y/3)
\]
An animation illustrating the steps in the method above, in the construction of the \emph{Blugold Fireworks} design,
can be found at the link given in \cite{ref:BlugoldFireworksAnimation}.
\subsubsection{Rotational Symmetry and Tessellation of $\H$\label{subsubsec:RoSymAndTessellationOfH}}
\begin{figure}[!htb]\centering
\setlength{\unitlength}{1in}
\begin{picture}(4.5,2.1)
\put(0,0){\resizebox{4.5in}{!}{
\begin{tikzpicture}
\node[anchor=south west] (image) at (0,0) {\includegraphics[width=3in]{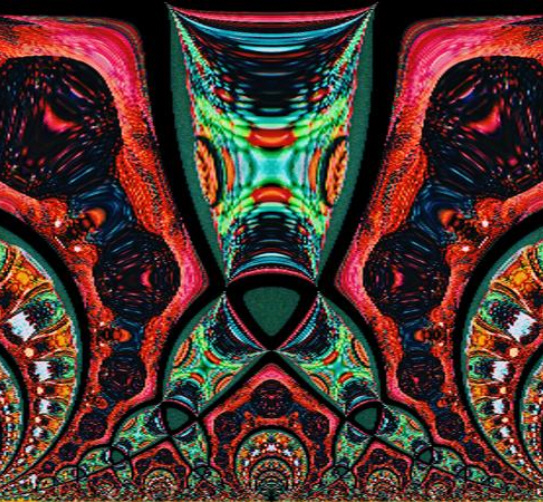}};
\begin{scope}[x={(image.south east)},y={(image.north west)}]
\draw[cyan, ultra thick] (0.5,0.55) -- (0.725,0.865); 
\draw[orange, ultra thick] (0.5,0.54) -- (0.5,0.4); 
\draw[orange, ultra thick] (0.5,0.4) -- (0.5,0.15); 
\draw[cyan, ultra thick] (0.24,0.28) arc (134:43:0.365);
\draw[cyan, dashed, ultra thick] (0.24,0.28) arc (134:180:0.365);
\draw[cyan, dashed, ultra thick] (0.24,0.28) arc (134:-.4:0.365);
\draw[yellow,ultra thick] (0.5,0.55) circle (1.16in); 
\fill[cyan] (0.5,0.55) circle (2.5pt); 
\fill[yellow] (0.5,0.38) circle (2.5pt); 
\fill[cyan] (0.5,0.15) circle (2.5pt); 
\draw[cyan, ultra thick] (1.5,0.55) -- (1.725,0.865); 
\draw[orange, ultra thick] (1.5,0.54) -- (1.5,0.4); 
\draw[orange, ultra thick] (1.5,0.4) -- (1.5,0.15); 
\draw[cyan, ultra thick] (1.23,0.28) arc (133:46:0.39);
\draw[cyan, dashed, ultra thick] (1.23,0.28) arc (135:180:0.39);
\draw[cyan, dashed, ultra thick] (1.23,0.28) arc (133:0:0.39);
\draw[yellow,ultra thick] (1.5,0.55) circle (1.16in); 
\fill[cyan] (1.5,0.55) circle (2.5pt); 
\fill[yellow] (1.5,0.383) circle (2.5pt); 
\fill[cyan] (1.5,0.15) circle (2.5pt); 
\node at (1.57,0.41) {\footnotesize$x\!+\! iy$}; 
\node at (1.365,0.55) {\footnotesize$x\!+\! iy\, \mathrm{cosh}\,\rho$}; 
\node at (1.47,0.27) {\footnotesize$\rho$}; 
\node at (1.35,0.32) {\footnotesize$\rho$}; 
\node at (1.62,0.32) {\footnotesize$\rho$}; 
\node at (1.6,0.13) {\footnotesize$x\!+\! iye^{-\rho}$}; 
\node[rotate=52] at (1.59,0.73) {\footnotesize$\mathrm{sinh}\,\rho$}; 
\end{scope}
\end{tikzpicture}}}
\end{picture}
\caption{\small Left: Illustration of circular region in $\H$ with center at point of $3$-fold rotational symmetry in $\H$. Right: Hyperbolic circle of radius $\rho$ centered at the point $x + iy\in \H$, and its Euclidean center and radius in $\EEL$.\label{fig:HyperbolicCircleIllustration}}
\end{figure}
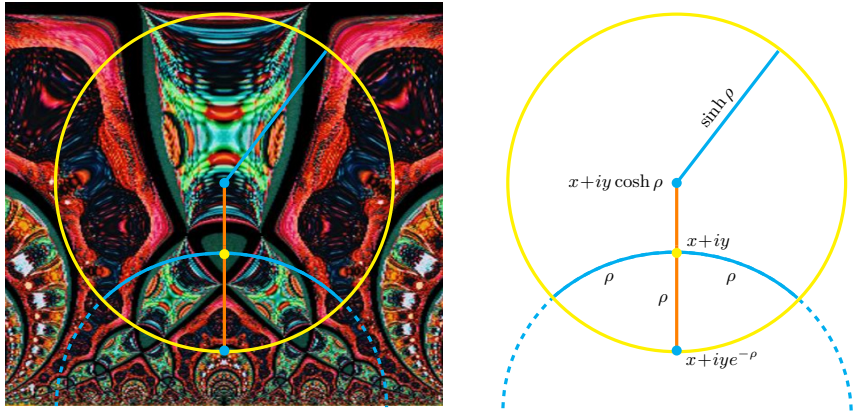

\begin{figure}[!htb]\centering
\setlength{\unitlength}{1in}
\begin{picture}(4.8,3.8)
\put(1.03,1.9){\resizebox{!}{1.9in}{
\begin{tikzpicture}
\fill[gray!20] (-0.97,5) -- (0.97,5) -- (0.97,2.03) -- (-0.97,2.03) -- cycle;
\fill[gray!20] (-0.97,2.04) -- (-0.245,2.04) -- (-0.97,1.85) -- cycle;
\fill[gray!20] (0.97,2.04) -- (0.245,2.04) -- (0.97,1.85) -- cycle;
\draw[dashed] (-4,0) -- (4,0);
\draw (-4,-0.1) -- (-4,0.1);
\draw (-2,-0.1) -- (-2,0.1);
\draw (0,-0.1) -- (0,0.1);
\draw (2,-0.1) -- (2,0.1);
\draw (4,-0.1) -- (4,0.1);
\node at (-3.9,-0.3) {\footnotesize$-2$};
\node at (-2,-0.3) {\footnotesize$-1$};
\node at (0,-0.3) {\footnotesize$0$};
\node at (2,-0.3) {\footnotesize$1$};
\node at (4,-0.3) {\footnotesize$2$};
\draw (-1,0) -- (-1,5);
\draw[orange, thick] (1,0) -- (1,5);
\draw (-3,0) -- (-3,5);
\draw (3,0) -- (3,5);
\draw (-4,2) arc (90:0:2);
\draw (-4,0) arc (180:0:2);
\draw (-2,0) arc (180:0:2);
\draw (0,0) arc (180:0:2);
\draw (2,0) arc (180:90:2);
\draw (-4,0) arc (180:0:0.67);
\draw (-3.33,0) arc (180:0:0.67);
\draw (-2,0) arc (180:0:0.67);
\draw (-1.33,0) arc (180:0:0.67);
\draw[cyan, thick] (0,0) arc (180:0:0.67);
\draw[cyan, thick] (0.67,0) arc (180:0:0.67);
\draw (2,0) arc (180:0:0.67);
\draw (2.67,0) arc (180:0:0.67);
\node at (-2,3.5) {\scriptsize$T^{-1}$};
\node at (-3.5,3.5) {\scriptsize$T^{-2}$};
\node at (0,3.5) {\scriptsize$\mathrm{Id}$};  
\node at (0,2.85) {\scriptsize$\F$};
\node at (2,3.5) {\scriptsize$T$};
\node at (3.5,3.5) {\scriptsize$T^2$};
\node at (-3.75,1.6) {\tiny$T^{-2}\mathcal{I}$};
\node at (-2,1.6) {\tiny$T^{-1}\mathcal{I}$};
\node at (0,1.6) {\tiny$\mathcal{I}$};
\node at (2,1.6) {\tiny$T\mathcal{I}$};
\node at (3.75,1.6) {\tiny$T^{2}\mathcal{I}$};
\node[rotate=55] at (-3.4,1.1) {\tiny$T^{-1}\mathcal{I}T\mathcal{I}$};
\node[rotate=-55] at (-2.5,1) {\tiny$T^{-1}\mathcal{I}T$};
\node at (-1.35,0.95) {\tiny$\mathcal{I}T\mathcal{I}$};
\node at (-0.65,0.95) {\tiny$\mathcal{I}T$};
\node at (0.65,1) {\tiny$\mathcal{I}T^{-1}$};
\node at (1.35,0.95) {\tiny$T\mathcal{I}T$};
\node[rotate=55] at (2.6,1.1) {\tiny$T\mathcal{I}T^{-1}$};
\node[rotate=-55] at (3.4,1.1) {\tiny$T^2\mathcal{I}T$};
\fill[yellow] (1,0.58) circle (2pt);
\fill[green] (1,1.73) circle (2pt);
\fill[red] (0,2) circle (2pt);
\end{tikzpicture}}}
\put(0,0){\resizebox{2.5in}{!}{
\begin{tikzpicture}
\node[anchor=south west] (image) at (0,0) {\includegraphics[width=3in]{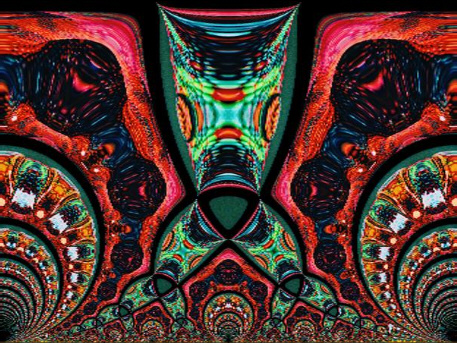}};
\begin{scope}[x={(image.south east)},y={(image.north west)}]
\draw[cyan, ultra thick] (0.33333,0.02) arc (180:0:0.325 and 0.41);
\draw[cyan, ultra thick] (0.66667,0.02) arc (0:180:0.325 and 0.41);
\draw[orange, ultra thick] (0.5,0.02) -- (0.5,0.975);
\fill[yellow] (0.5,0.38) circle (2.5pt);
\end{scope}
\end{tikzpicture}}}
\put(2.5,0){\resizebox{2.03in}{!}{ 
\begin{tikzpicture}
\node[anchor=south west] (image) at (0,0) {\includegraphics[width=3in]{reptilefortikz.jpg}};
\begin{scope}[x={(image.south east)},y={(image.north west)}]
\fill[yellow] (0.501,0.378) circle (3pt);
\fill[yellow] (0.328,0.1718) circle (1.5pt);
\fill[yellow] (0.2235,0.1025) circle (1pt);
\fill[yellow] (0.6748,0.17) circle (1.5pt);
\fill[yellow] (0.777,0.1025) circle (1pt);
\fill[red] (0.501,0.66) circle (3pt);
\fill[red] (0.38,0.275) circle (2.5pt);
\fill[red] (0.2625,0.148) circle (2pt);
\fill[red] (0.618,0.2758) circle (2.5pt);
\fill[red] (0.7375,0.15) circle (2pt);
\draw[cyan,very thick] (0.29,0.98) arc (180:263:0.78 and 0.965);
\draw[cyan,very thick] (0.71,0.98) arc (0:-83:0.78 and 0.965);
\end{scope}
\end{tikzpicture}}}
\end{picture}
\caption{\small Top: Tessellation of $\H$ generated by $T$ and $\I$.
Regions are labeled by transformations that produce them from shaded region $\mathcal{F}$.  Bottom Left:
Location of center of $3$-fold symmetry (yellow dot). This center lies at the intersection of three geodesics, shown in blue and orange in this image and in the tessellation above it. Bottom Right: Centers of $3$-fold (yellow dots) and $2$-fold (red dots) symmetries located along horocycles.\label{fig:DedekindTesselationIllustration}}
\end{figure}

The designs we have created relate to a number of other additional aspects of the geometry of $\H$.
For example, on the left of \Fig~\ref{fig:HyperbolicCircleIllustration} we show a clip from the
third design in \Fig~\ref{fig:HyperbolicDesigns}. Overlayed on this clip is a circle surrounding
a point of $3$-fold rotational symmetry in $\H$. It is important to note that this $3$-fold rotational\
symmetry exists within $\H$ and not within $\EEL$. To see that we do have $3$-fold symmetry in $\H$, in contrast
to what we are used to seeing with Euclidean geometry, we need to discuss some further ideas from hyperbolic
geometry. We first begin by describing the significance of the yellow circle on the left of \Fig~\ref{fig:HyperbolicCircleIllustration}. Its center in $\EEL$ is marked by a blue dot. This circle is the locus
of points that are a fixed distance $\rho$ from the yellow dot in $\H$. To be precise, we have
the following theorem.
\begin{theorem}\label{thm:HyperbolicCircles}%
\textsl{For fixed $\rho > 0$, the locus of points $C_\rho (x + iy)$ that are distance $\rho$ in $\H$ from a fixed
point $x + iy$ is equal to the circle in $\EEL$ with center $x + iy \cosh \rho$ and radius
$\sinh \rho$.} {\small [See \Fig~\ref{fig:HyperbolicCircleIllustration}.]}
\end{theorem}
\begin{small}\begin{proof}
Begin by supposing that the fixed point is $i$. Using the distance formula $d(i y_1, i y_2) = |\ln(y_2/y_1)|$,
the two points $i e^\rho$ and $i e^{-\rho}$ are both distance $\rho$ from $i$ in $\H$. So,
$i e^\rho, i e^{-\rho} \in C_\rho (i)$. The midpoint
on the $i$-axis in $\EEL$ between these two points is $i\cosh \rho$, and it is Euclidean distance $\sinh \rho$
from both points.  Now, map $C_\rho (i)$ in $\H$ to $C_\rho (0)$ in the unit disc $\D$, using the isometry $f(z) = (iz + 1)\, /\, (z + i)$ from $\H$ to $\D$ described
below in Section~\ref{subsubsec:MappingsToDiskD}. As discussed in that same section, the metric differential
$ds_{\hbox{\tiny $\D$}}$ has rotational invariance about $0$, and therefore
$C_\rho(0)$ is a Euclidean circle about $0$ (although its radius is \Emph{not} equal to $\rho$). Then map $C_\rho(0)$ back
to $C_\rho (i)$ in $\H$ using the isometry $f^{-1}(z) = (z + i)\, /\, (iz + i)$ from $\D$ to $\H$. Because $f^{-1}(z)$ is a linear fractional transformation, it maps Euclidean circles to Euclidean circles, hence $C_\rho(i)$ is a circle in $\EEL$. Since we found its center and radius must be $i\cosh \rho$ and $\sinh \rho$, we have proved the result for fixed point $i$. Conjugating with $\T_x \M_y$, we get the result for fixed point $x + iy$.
\end{proof}\end{small}
The yellow dot, located within a triangular region with curved edges on the left of \Fig~\ref{fig:HyperbolicCircleIllustration}, is a center for a $3$-fold rotation in $\H$. It is located at
$1/2 + (\sqrt{3}/6)i$. Before we show that it is a center for a $3$-fold rotation in $\H$, we show that the
point $z_1 = 1/2 + (\sqrt{3}/2)i$ is also a center for a $3$-fold rotation in $\H$. This point $z_1$ is marked
by a green dot on the graph on the top of \Fig~\ref{fig:DedekindTesselationIllustration}. It is related to
a tessellation of $\H$ that we will discuss shortly. For now, observe that it lies at a vertex of a region
labeled $T\I$ in the tessellation. We will show that $z_1$ is a fixed point for $T\I$, and that $T\I$ has order
$3$ in the group $\MGZ$. We have
\[
T\I (z_1) = \frac{-1}{1/2 + (\sqrt{3}/2)i} + 1\; = 1/2 + (\sqrt{3}/2)i
\]
so $z_1$ is a fixed point for $T\I$. Moreover, we have
\[
T\I (z) = \frac{z-1}{z}\, ,\qquad \left(T\I\right)^2 (z) = \frac{-1}{z-1}\, ,\qquad \left(T\I\right)^3 (z) = z
\]
which shows that $T\I$ has order $3$ in $\MGZ$.

We now turn to $z_2 = 1/2 + (\sqrt{3}/6)i$. We observe that $z_2 = \I T^{-2} z_1$. But, $\I T^{-2}$ is an isometry.
Hence we can apply the following Lemma:
\begin{lemma}\label{lem:ConjugatingIsometriesAndFixedPoints}\textsl{%
Suppose $z$ is a fixed point for $g\in\S_\H$, and that $g$ has finite order $k$. If $f \in\S_\H$, then
$f\circ g \circ f^{-1}$ has fixed point $f(z)$ and order $k$.}
\end{lemma}
\begin{small}\begin{proof}
We find that $(f\circ g \circ f^{-1}) \circ f(z) = (f\circ g)(z) = f(z)$ so $f(z)$ is a fixed point. Moreover,
$(f\circ g \circ f^{-1})^j = f\circ g^j \circ f^{-1}$ for any integer $j\geq 0$. When $j = k$, we
have $g^k = \mathrm{Id}$, so $(f\circ g \circ f^{-1})^k = f\circ f^{-1} = \mathrm{Id}$. Also, when
$j < k$, if $f\circ g^j \circ f^{-1} = \mathrm{Id}$, then we would have $g^j = f^{-1}\circ f \, = \mathrm{Id}$
and that would contradict $k$ being the order of $g$. Consequently, $f\circ g \circ f^{-1}$ has order $k$.
\end{proof}\end{small}\noindent
Applying the Lemma, we see that $1/2 + (\sqrt{3}/6)i = \I T^{-2} z_1$ is a fixed point for
$h = \I T^{-2} \circ T\I \circ (\I T^{-2})^{-1}$ and $h$ has order $3$.

Returning to the left of \Fig~\ref{fig:HyperbolicCircleIllustration}, the significance of the yellow dot and the circle
enclosing it can now be explained in terms of rotation in $\H$. Since the yellow dot corresponds to the
fixed point $z_2$ for the isometry $h$ of order $3$, it is analogous to a $3$-fold rotation in $\EEL$. In fact,
for points sufficiently close to $z_2$, the metric differential $ds = ds_{\Esub}/y$ is approximately equal to
a multiple of $ds_{\Esub}$. Consequently, the isometry $h$ is acting like a $3$-fold rotation in $\EEL$ in the
limit of approaching $z_2$. The enlargement of parts of the design as one rotates towards the vertical corresponds
to what we see in $\EEL$. By Theorem~\ref{thm:IsometriesInH} we know that area is preserved by isometries. Hence,
in $\H$, the upper arm of the figure along the vertical direction, has exactly the same area as each of the
two lower arms of the figure (extending out from the edges of the curved triangle).

There are also many other centers of $3$-fold rotation in $\H$ that are illustrated in this design. If we
let $z_k = \I T^k z_1$ for $k = -3, -4, -5, \dots$, then Lemma~\ref{lem:ConjugatingIsometriesAndFixedPoints}
implies that we have centers of $3$-fold rotations at each $z_k$. Since each point $T^k z_1 = z_1 + k$ lies
along a horizontal horocycle in $\EEL$, applying the inversion $\I$ maps them to a circular horocycle
in $\EEL$. We can see some of these points $z_k$ in the image on the bottom right of \Fig~\ref{fig:DedekindTesselationIllustration}.
They are lying above the blue curve, extending downwards towards the bottom left corner, which is a slightly lower
horocycle belonging to the same family of horocycles tangent to $0$ in $\R$. Finally, let $z_k = T\I T^k z_1$ for
$k = 2, 3, 4, \dots$. This produces another collection of centers of $3$-fold rotations that are on a
second horocycle that moves away to the right of $z_1$ towards $1$ in $\R$. Some of these centers are
visible in the figure as well. For all of these centers, we observe repetitions of the curved
triangle containing $z_1$ but at smaller scale. They are only smaller scale in $\EEL$. In $\H$, isometries
preserve area, so these curved triangles all have the same area in $\H$ as the one containing $z_1$, i.e.,
they are of the same scale in $\H$.

In addition to $3$-fold centers, there are numerous centers of $2$-fold rotational symmetries in $\H$. The
point $z = i$ is a fixed point for $\I$ which has order $2$. It is shown as a red dot in the tessellation
at the top of \Fig~\ref{fig:DedekindTesselationIllustration}. Applying isometries to $i$ we obtain
sequences of $2$-fold centers lying along horocycles. These centers of symmetry are plotted as red
dots on the image shown on the right of \Fig~\ref{fig:DedekindTesselationIllustration}.

We have discussed the relation between our design and the tessellation shown at the top of \Fig~\ref{fig:DedekindTesselationIllustration}.
We shall now discuss this tessellation in more detail. Equation~\eqref{eq:DecompositionWithTandI} shows how any $f\in\MGZ$
can be written in terms of powers of $T$ interlaced with $\I$.  Starting from a
\emph{fundamental domain}, indicated by the shaded region $\F$ in the figure, and applying compositions
of $\I$ with powers of $T$ generates this tessellation of $\H$. In fact, the construction of the tessellation reproduces
all the possible combinations of powers of $T$ interlaced with $\I$ in Equation~\eqref{eq:DecompositionWithTandI}.
We now make all of these ideas precise with the following theorem.
\begin{theorem}\label{thm:TesselationOfH}\textsl{%
The modular group $\MGZ$ generates a tessellation of $\H$ via
\begin{equation}\label{eq:TesselationOfH}
\H = \bigcup_{f\in\MGZ} f\bigl(\F\bigr)
\end{equation}
where $\mathcal{F} = \{z\in\H\;\colon\, |z| \geq 1, |\Re{z}| \leq 1/2\}$, and every pair of
regions $f\bigl(\F\bigr)$ and $g\bigl(\F\bigr)$ for $f\neq g$ have disjoint interiors.}
\end{theorem}
\begin{small}\begin{proof}
First, we show that $\H = \bigcup_{f\in\MGZ} f\bigl(\F\bigr)$. Let $z$ be an arbitrary element in $\H$.
We will show that there is a $g(z) = \hbox{\normalsize $\frac{jz + k}{mz+n}$}\,\in\MGZ$ for which $w = g(z) \in \F$. Then
we will have $f(w) = z$ for $f = g^{-1}$, which will establish the decomposition of $\F$ in Equation~\eqref{eq:TesselationOfH}.
Our main tool is
\begin{equation}\label{eq:ImagConstraint}
\Im{g(z)} = \Im{z}\, / \, |mz + n|^{2}
\end{equation}
which was shown for all isometries at the beginning of the proof of Theorem~\ref{thm:IsometriesInH}.
Since there are only finitely many $m, n \in\Z$ for which $|mz + n| \leq 1$, 
it follows from \eqref{eq:ImagConstraint} that there are only finitely many $g \in\MGZ$ for which
$\Im{g(z)} \geq \Im{z}$. Therefore, we can choose a $g\in\MGZ$ for which
$\Im{g(z)}$ is maximal. If $|\Re{g(z)}| > 1/2$, then we may compose
$g$ with some power of $T$ so that $|\Re{g(z)}| \leq 1/2$. Therefore, without loss
of generality, we assume that $|\Re{g(z)}| \leq 1/2$. Then we must have
$|g(z)| \geq 1$, because if $|g(z)| < 1$ we would have $0 < \Im{g(z)} < 1$, hence
$\Im{\bigl[\I g(z)\bigr]} = \mathrm{Im}\bigl[\hbox{\normalsize $\frac{-1}{g(z)}$}\bigr] > \Im{g(z)}$ and
that contradicts the maximality of $\Im{g(z)}$. Thus, $|\Re{g(z)}| \leq 1/2$
and $|g(z)| \geq 1$, and so $w = g(z) \in \F$.

Second, we prove that $f\bigl(\F\bigr)$ and $g\bigl(\F\bigr)$ have disjoint interiors when $f\neq g$. But this
is equivalent to $\F$ and $g\bigl(\F\bigr)$ having disjoint interiors for every
$g(z) = \hbox{\normalsize $\frac{jz + k}{mz+n}$} \in\MGZ$ not equal to $\mathrm{Id}$. Suppose first that $m = 0$.
Then $g(z) = z \pm k$, for $k\neq 0$, hence $\F$ and $g\bigl(\F\bigr)$ have disjoint interiors. Now, suppose
$m\neq 0$. Let $z$ be in the interior of $\F$. Hence $|\Re{z}| < 1/2$ and $|z| > 1$. We then have
\begin{align*}
|mz + n|^2 &= m^2 |z|^2 + 2(\Re{z})mn + n^2\\
&> m^2 + n^2 - |mn| \; = (|m| - |n|)^2 + |mn|
\end{align*}
Since $m \neq 0$, the strict lower bound $(|m| - |n|)^2 + |mn|$ is a positive integer. Consequently,
$|mz + n|^2 > 1$ and so $\Im{g(z)} < \Im{z}$. Since
$g^{-1}(w) = \hbox{\normalsize $\frac{nw - k}{-m w + j}$}$, if $w = g(z)$ were in the interior of
$\F$, then the same argument yields $\Im{g^{-1}(w)} < \Im{w}$. Hence
$\Im{z} < \Im{g(z)}$, and this contradiction shows that the interiors of $\F$ and $g\bigl(\F\bigr)$ are disjoint.
\end{proof}\end{small}\noindent
A remarkable feature of the animation in~\cite{ref:BlugoldFireworksAnimation} is how new features are added
near the bottom of the screen that clearly correspond to the bottom portions of the tessellation shown
in \Fig~\ref{fig:DedekindTesselationIllustration}.
\subsubsection{Mappings to the Disk $\D$\label{subsubsec:MappingsToDiskD}}
We have emphasized the hyperbolic upper half-plane as a model for non-Euclidean geometry. Another model, which
is equivalent to $\H$, uses the unit disc $\D = \{z\in\C\;\colon\;|z| < 1\}$ as underlying set of points.
Needham~\cite[pp.~317--318]{Needham} shows that the map $f(z) = \hbox{\Large $\frac{iz + 1}{z + i}$}$ is a conformal map
from $\H$ to $\D$, and it induces a metric differential $ds_{\hbox{\tiny $\D$}}$ on $\D$ given by
\begin{equation}\label{eq:DiskMetricDifferential}
ds^2_{\hbox{\tiny $\D$}} = \frac{2}{1 - |z|^2}\, ds^2_{\hbox{\tiny $\mathbb{E}$}}
\end{equation}
where $ds^2_{\hbox{\tiny $\mathbb{E}$}} = dx^2 + dy^2$ is the Euclidean metric differential for $\D$ as a subset of the
Euclidean plane $\mathbb{E}$. With this metric for $\D$, the map $f\colon \H \to \D$ is an isometry. Consequently, it maps geodesics in $\H$ to geodesics in $\D$.
The geodesics in $\D$ lie on either arcs of circles that intersect the unit circle at right angles, or diameters in $\D$ (which also intersect the unit circle at right angles). An excellent treatment of this disc model for non-Euclidean geometry can be found in Krantz~\cite{ref:Krantz}. More details about the geometry, including tessellations of $\D$, are in
Climenhaga and Kotek~\cite{ref:ClimKotek}.

\begin{figure}[!b]\centering
\setlength{\unitlength}{1in}
\begin{picture}(6.5,2.95)
\put(0,0.2){\includegraphics[width=2.7in]{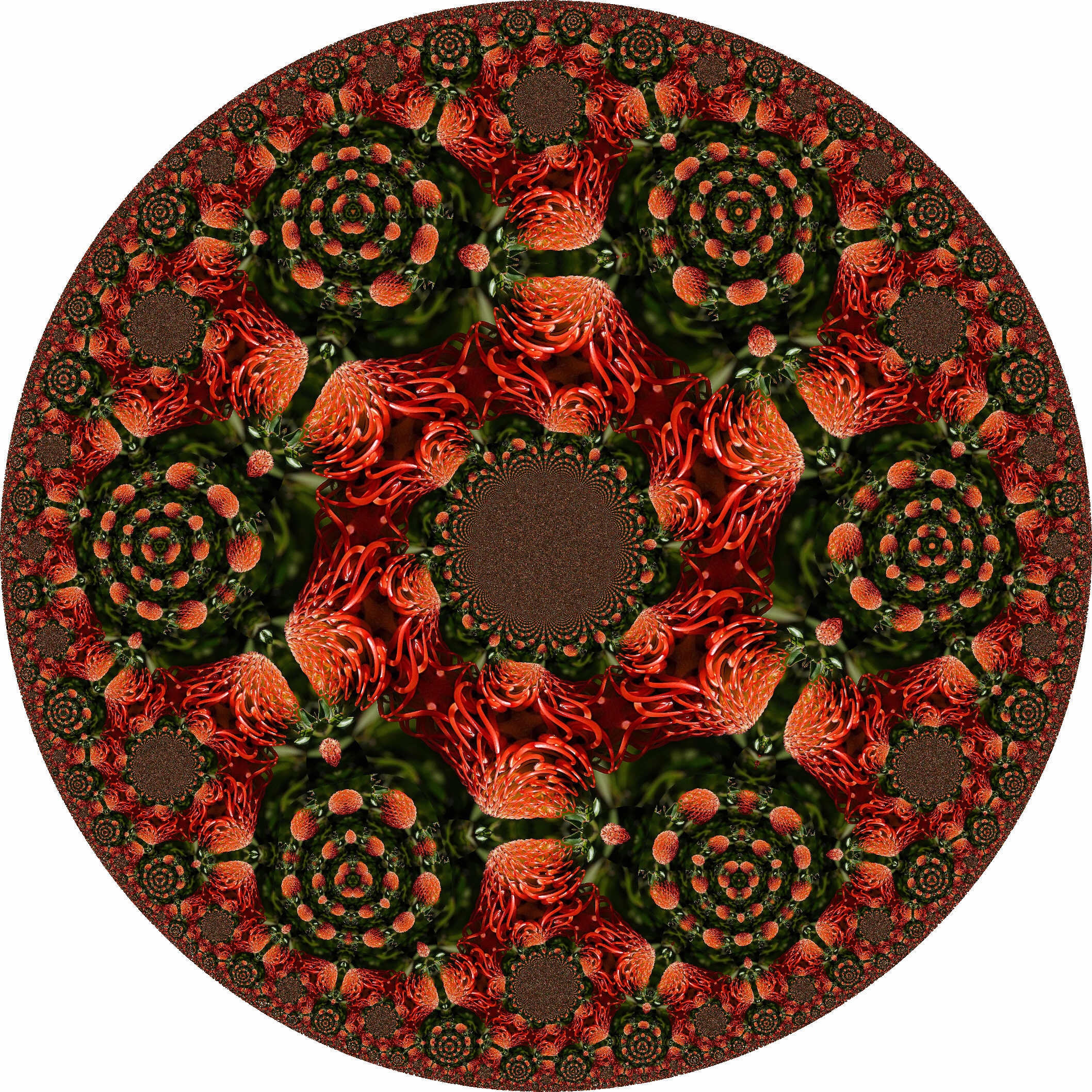}} 
\put(2.8,0.2){\includegraphics[width=2.7in]{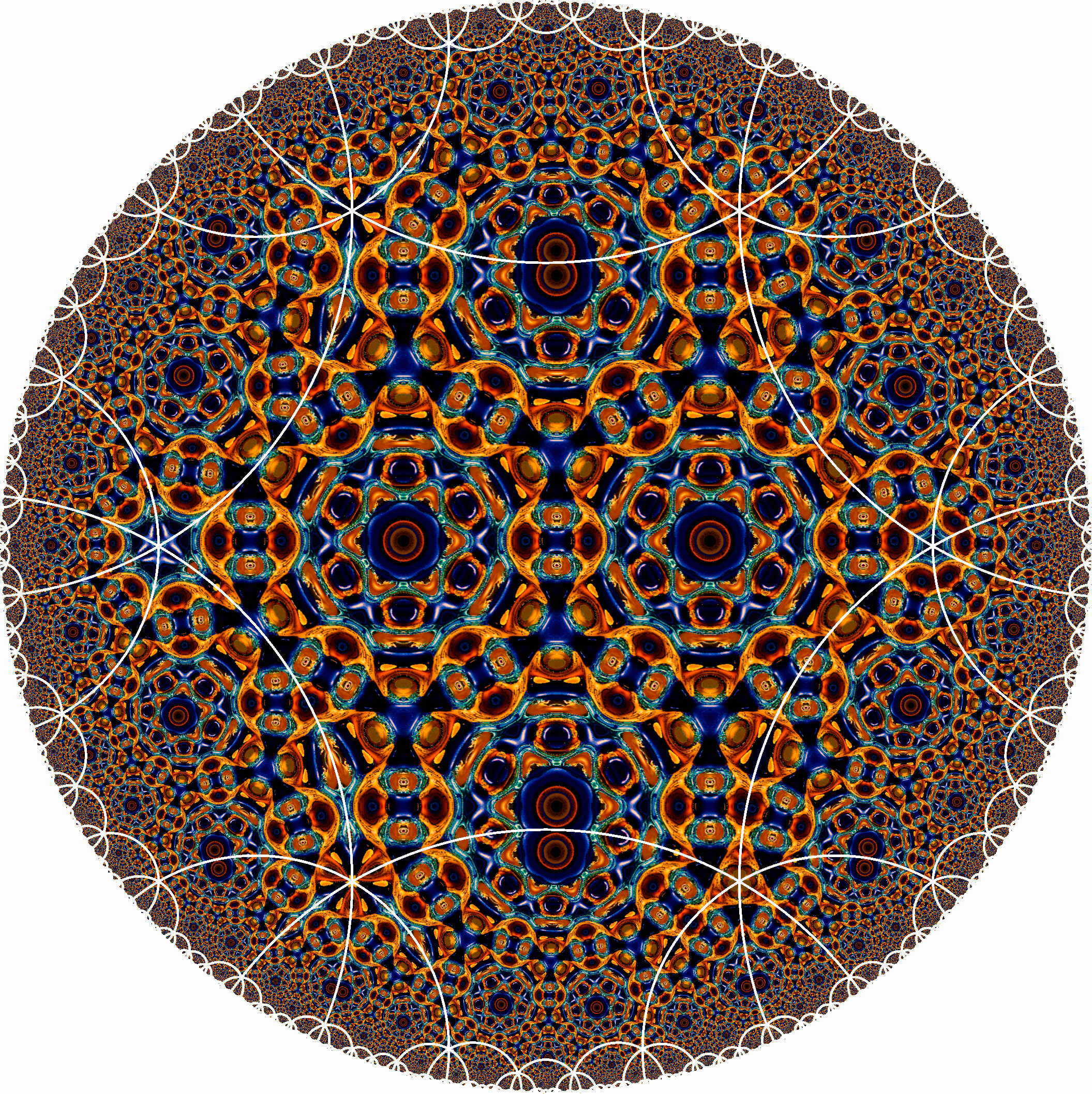}} 
\put(5.6,1.95){\includegraphics[width=.85in]{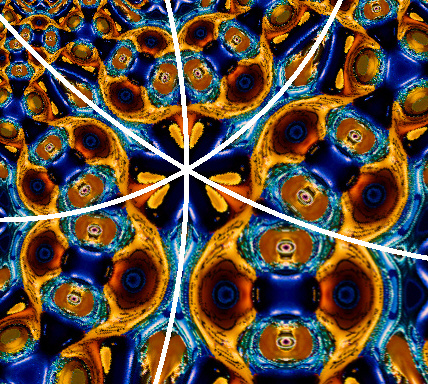}}
\put(5.6,.9){\includegraphics[width=.85in]{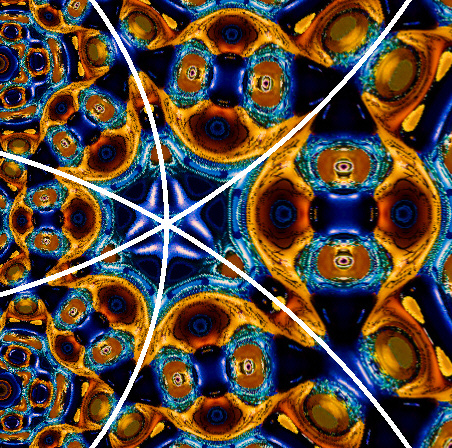}}
\put(1.1,0){\scriptsize \textsf{Design~1}}
\put(4,0){\scriptsize \textsf{Design~2}}
\end{picture}
\caption{\small Two designs created with \cite{ref:DiskDesignApp}. {\footnotesize\textsf{Design~2}} includes
tessellating curves, shown in white.  To the right of {\footnotesize\textsf{Design~2}} are two zooms showing $3$-fold symmetries in
the hyperbolic disc $\D$.
\label{fig:DiscMappings}}
\end{figure}

There is an abundance of artistic designs in $\D$ that have already been created. The designs by Escher, using
curved polygonal tessellations of $\D$ are surely the most famous~\cite{Escher}. The mathematics of creating curved polygonal tessellations of $\D$ was worked out by Coxeter~\cite{ref:Coxeter}. An entertaining app for creating your own designs, using Coxeter's tessellations of $\D$, can be found at the link in~\cite{ref:DiskDesignApp}. In \Fig~\ref{fig:DiscMappings} we show two designs we created with this app.
\noindent
{\small\textsf{Design~1}} used one of our \emph{Waratah} flower rosettes as source image. It clearly retains the original $6$-fold rotational symmetry of the original rosette. More importantly, from an artistic standpoint, it contains an ambiguity of \emph{figure-ground} relations. If you view {\small\textsf{Design~1}} from a far distance, the design features a star with six arms emanating from the central region of the design. However, when you move up close, this star fades into the background and six green/red flower-like regions closer to the circular boundary of the disc are more prominent.

{\small\textsf{Design~2}} used one of our $6$-fold symmetric wallpaper designs as source image.
This latter design is particularly interesting in that it includes curved polygons (shown in white) that form a
tessellation of $\D$. The app~\cite{ref:DiskDesignApp} created this design by loading a clip of the source image
into a central curved hexagon, shown in white at the center of {\small\textsf{Design~2}}. During the design process,
it successively displays the iteration of isometries of this central, fundamental, hexagonal region to fill out
the rest of the curved hexagons that are tessellating the disc. As with the tessellation of $\H$ we discussed above,
there are centers of $2$-fold and $3$-fold symmetries in the completed design.
It is easy to spot centers for $2$-fold and $3$-fold \emph{hyperbolic} rotational symmetry near the top and bottom of the image. The $3$-fold centers are located at intersections of the tessellating curves (see the two zooms on the right of \Fig~\ref{fig:DiscMappings}), just as they occur at such intersections for the tessellation of $\H$. Notice that
the original $6$-fold symmetry of the source image is lost in {\small\textsf{Design~2}}, due to the clipping of only
a part of the source image within the central curved hexagon. From an artistic standpoint, we like this \emph{symmetry
breakage} in {\small\textsf{Design~2}}.

It is also an important fact that both $1/(1 - |z|^2)$ and $ds^2_{\hbox{\tiny $\mathbb{E}$}} = dx^2 + dy^2$ are invariant under rotations centered at $0$. Since $|e^{i\theta} z|^2 = 1$, we have the invariance of $1/(1 - |z|^2)$ by rotation by $\theta$, and the invariance of $ds^2_{\hbox{\tiny $\mathbb{E}$}}$ under rotation holds because rotation by $\theta$ can be expressed as an orthogonal matrix. Since $ds^2_{\hbox{\tiny $\D$}}$ is the product of $1/(1 - |z|^2)$ and $ds^2_{\hbox{\tiny $\mathbb{E}$}}$, it
is invariant under rotation about the origin. This rotational invariance is exhibited near the centers of both {\small\textsf{Design~1}} and {\small\textsf{Design~2}} in \Fig~\ref{fig:DiscMappings}, where these designs retain the rotational symmetry about the origin enjoyed by the rosette and wallpaper designs used as their source images. Similar $3$-fold and $2$-fold hyperbolic rotational symmetries, and $4$-fold Euclidean rotational symmetry about the origin, are even more clearly evident in Escher's \emph{Circle Limit III} woodcut~\cite{Escher} (if one ignores the different colors of the fishes).

\begin{figure}[!b]\centering
\setlength{\unitlength}{1in}
\begin{picture}(2.1,2.15)
\put(0,0){\includegraphics[width=2.1in]{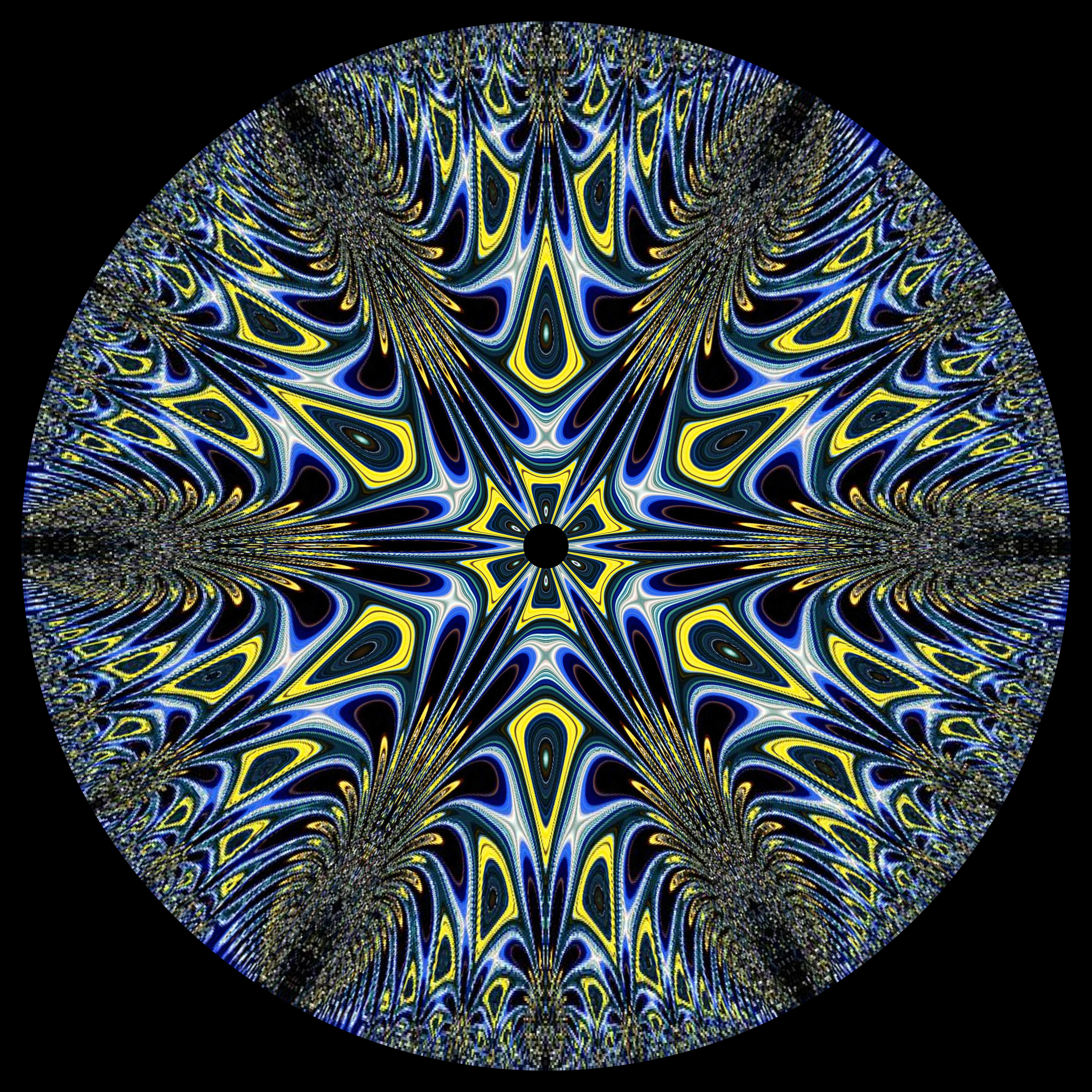}} 
\end{picture}
\caption{\small Design created using $\log z + \log(-z)$. 
\label{fig:DiscMappingLogProcess}}
\end{figure}

We have also produced images using conformal maps from $\H$ to $\D$. For example, Stein~\cite[p.~211]{ref:Stein} shows that
$\log z$ defines a conformal map from the half-disc $\{z = x + iy \;\colon\, |z| < 1, y > 0\} \subset \D$
to the half-strip $\{w = u + iv\;\colon\, u < 0, 0 < v < \pi\} \subset \H$.
With this idea in mind, we created a design in $\D$ using a multiple of $f(z) = \log z \, + \, \log(-z)$ and
\emph{Blugold Fireworks} as color map.
It is shown in \Fig~\ref{fig:DiscMappingLogProcess}. This design clearly shows $6$-fold rotational symmetry about
the origin in $\D$ and reflection symmetries through six diameter geodesics.
We intend to continue exploring mappings from $\H$ to $\D$, including ones that more fully exploit the isometric equivalence
of $\H$ and $\D$.
\noindent
\subsection*{Conclusion}\label{subsec:ConclusionHypDesigns}
We have shown a number of symmetric designs generated by the application of complex analysis to the
non-Euclidean geometry of the hyperbolic upper-half plane and the hyperbolic disc. The history of the mathematics of
hyperbolic geometry is a fascinating one. References for this history include Needham~\cite{Needham} and
Greenberg~\cite{ref:Greenberg}. Penrose~\cite{ref:Penrose}[Chap.~2] has some fascinating insights. Mathematics related to hyperbolic geometry continues right up to the present day, see~e.g., Adams~\cite{ref:Adams}.

\end{document}